\documentclass[10pt,final,twoside]{amsart}

\usepackage{booktabs,multirow,array,colortbl,tabularx}
\usepackage{flafter,pdflscape}
\usepackage{soul,cite,lineno,cancel,verbatim,ulem}
\usepackage{slashbox,xcolor,relsize,tikz,tikz-cd,rotating,lscape,adjustbox}
\usepackage{latexsym,amsfonts,amsthm,mathtools}
\usepackage{enumitem}
\usepackage{graphicx}
\graphicspath{ {images/} }
\evensidemargin=-
.5cm
\newcommand{\xb}{\mathbf{x}}

\newcommand{\xg}{x_{\star}}  
\newcommand{\xs}{x^{\star}}

\newcommand{\yb}{\mathbf{y}}
\newcommand{\yg}{y_{\star}}  
\newcommand{\ys}{y^{\star}}

\newcommand{\specialcell}[2][c]{%
	\begin{tabular}[#1]{@{}c@{}}#2\end{tabular}}

\newtheorem{theorem}{Theorem}[section]

\newtheorem{remark}[theorem]{Remark}

\begin{document}
	\title{
A comparison among a fuzzy algorithm for image rescaling  with other methods of digital image processing
} 
	\author{Danilo Costarelli}
	\author{Anna Rita Sambucini*}
\address{
	Danilo Costarelli\newline
	\indent  University of Perugia \newline
	\indent Department of Mathematics and Computer Sciences\newline
	\indent 1 Via Vanvitelli, 06124, Perugia, Italy\newline
	\indent ORCID: {0000-0001-8834-8877}}
\email{danilo.costarelli@unipg.it}

\address{
	Anna Rita Sambucini\newline
	\indent  University of Perugia \newline
	\indent Department of Mathematics and Computer Sciences\newline
	\indent 1 Via Vanvitelli, 06124, Perugia, Italy\newline
	\indent ORCID: {0000-0003-0161-8729}}
\email{anna.sambucini@unipg.it}


	\begin{abstract}
The aim of this paper is to  compare  the fuzzy-type algorithm for image rescaling introduced by Jurio et al., 2011, 
quoted in the list of references, with some other existing algorithms such as  the classical bicubic algorithm and 
the  sampling Kantorovich (SK) one. 
Note that  the SK algorithm is a recent tool for image rescaling and enhancement that  has been
revealed to be useful in several applications to real world problems, while  the  bicubic algorithm is widely known in the literature.
\\
 A comparison among the  abovementioned algorithms (all implemented   in the   MatLab programming language)
 was performed in terms of suitable similarity  indices
 such as the Peak-Signal-to-Noise-Ratio (PSNR)  and the likelihood index $S$. 
		\vspace{5pt}
		
		\vspace{5pt}
	\end{abstract}
	
\maketitle
\pagestyle{myheadings}
\markboth{D. Costarelli, A.R. Sambucini}{A comparison among a fuzzy algorithm for image rescaling ...}


\section{Introduction}\label{s-intro}
Images are  indispensable  tools in concrete life, as well as in  various 
 fields of research and they   have a concrete impact on daily life.
The most common scientific applications of image processing are in medicine, in
which some instrumental tests such as CT  and  MRI,  are  helpful  for the diagnosis
of various diseases, remote sensing, in which the use of satellite images allows
the study of phenomena (climatic, tectonic, etc) linked to natural events; astronomy, biology and  many other fields.

In real world applications  digital images are essential tools for studying concrete problems 
since  they  provide visual and numerical representations of an observation
or a measurement. Namely, they constitute 
a synthesis of information concerning one or more characteristics of the 
problem under consideration.

The acquisition of a digital image from a camera or a diagnostic device is a physical process that allows  
  the conversion of  measured data  into  two or three dimensional discrete signals/images. 

 During this phase, the acquisition tools, which are obviously endowed with their own sensitivity 
and by their own procedure of data conversion, allow     the reconstruction of
 a digital image that is obviously characterized by a natural degree of approximation and therefore of uncertainty, i.e., 
it is not always possible to establish the gray levels of a region of pixels perfectly or to precisely detect 
 geometric shapes   characteristics, such as edges of particular interest. 

These facts can be translated into the constructionof a matrix of pixels in which the value of each element represents a 
"good approximation” of the real gray level (luminance) in the gray scale.

When situations of this type are present, it is possible to use fuzzy set theory to represent and elaborate vague and imprecise
 concepts and apply  a fuzzy algorithm for digital image processing, as  in, 
 for example  \cite{BBPF2009,Mesiar,latorre,danilo,xdanilo}. 
 Moreover, we also know that fuzzy theory is a fundamental tool for several topics, such as probability
 (see, e.g., \cite{ADT1,DAT1,STV1}) and many others,
hence it is not surprising to find a close connection  between
 digital images and their processing.
 Recently  multifunctions have also  been applied 
for convergence results in this setting see also \cite{dipiazza1,dipiazza2,dipiazza3}.

On the other hand, it is also known that any image is a multivariate discontinuous signal, where the possibility  of visualizing
 the contours and edges in the figures is due to the presence of meaningful jumps of gray levels in the grayscale;
 this is the motivation why Signal Theory has been successfully applied to process digital images. 
 Indeed, in the last ten years, several models for concrete applications in the field of medicine and engineering have been developed thanks 
to the use of the  SK algorithm ( e.g., \cite{CLCOGUVI1,Poland1}). 
The main purpose of the SK algorithm is 
 to rescale images,  by acting as low-pass filter and hence contrasting the appearance of  noise. 
The SK algorithm is the  numerically optimized implementation of the sampling Kantorovich operators
 (from which  the acronym SK),  widely studied in Approximation Theory, since it is very suitable 
 for reconstructing non-necessarily continuous signals (hence images, see \cite{CSV3,COVI}).  
 The aim of this paper is to compare the fuzzy-type algorithm introduced in \cite{Mesiar} with the classical bicubic interpolation method,
 widely used in the literature  e.g., \cite{LIU} and the
 above described SK algorithm. 
The above algorithms  were implemented  in the  MatLab programming language, and 
 comparisons were performed by means of several numerical 
 tests performed on a suitable dataset of images of different types.

 To quantitatively evaluate the results, we introduced two similarity  indices known in the literature.
We considered the Peak-Signal-to-Noise-Ratio (PSNR) \cite{TAN}, and the likelihood index S considered in \cite{BPB2007}.
 Finally,   a comparison in   terms of CPU time employed by the three considered algorithms 
 was also carried out  for the best approximations.


\section{The Interval valued Fuzzy point of view}\label{due}
A grayscale digital image of dimensions $n \times m$ (i.e., with $n$ rows and $m$ columns) is a matrix $Q$ of dimensions $n \times m$,
 where the element of position $(i,j)$ in the matrix, denoted by $q_{i,j}$, represents the intensity of the pixel in the gray scale (luminance).
We observe that it is not restrictive to work only with grayscale images, since operating on a color image is  similar to  doing  so on 3 grayscale images. 
 For colour images three matrices are used which, for each pixel, assume integer values in the range $[0,255]$ 
with respect to the red, green and blue colours (RGB channels, see, for example, \cite{GON}).

  The  luminance value  $q_{i,j}$ at  point $(i,j)$  are normalized  to obtain values in the range $[0,1]$.
To simplify the notation, we will always indicate them with the same symbol.\\
In \cite{Mesiar},  Jurio,  Paternain,  Lopez-Molina,  Bustince,  Mesiar and  Beliakov proposed  a model 
 associated to a grayscale image and an interval valued fuzzy set  to construct  a magnification algorithm that 
 considers  the luminance values in a neighbourhood of  each pixel of the image. 

The type of operator they  use is of spatial type, namely, to determine the value of the destination pixels,
 not only the value of the pixel in the original image  but also the value of some
pixels close to it (in a neighbourhood of it)  will be considered.

The key idea of this rescaling algorithm  (proposed by Jurio et al.) is to associate  an interval membership to each pixel. 
The  parameter $\delta$ is fixed a priori; when $\delta$ increases the length of the interval increases, so more values of  the intensities
 of the pixels  close to the assigned intensity are  considered.
In this way, a new block is constructed for each pixel
of the image, and the central pixel of the block  maintains
the luminance of the original pixel. To fill the rest of the pixels in the  newly generated block,
the relationship between the luminance of the pixel in the original image 
and that of the pixels "near" to  the pixel was used.\\

 To define the interval-valued membership of $q_{i,j}$ 
 let $ L([0,1])$ be the family of all closed intervals in $[0,1]$, namely 
  \[ L([0,1]) := \{ \xb =[\xg,\xs] : \, (\xg, \xs) \in [0,1]^2\,\,\, \& \,\,\, \xg \leq \xs\,   \},\] 
 with the following partial order relation:
 $\xb \leq_L \yb$ if   $\xg  \leq \yg$ and $\xs \leq \ys$ (this is a lattice order between closed intervals; see, for example, \cite{Mesiar}).   
For every closed interval $\xb:=[\xg,\xs]$ in $L([0,1])$, let  $W(\xb) := \xs - \xg$ be its length.\\

Therefore, an  interval-valued membership of $q_{i,j}$  is an  
 interval valued fuzzy set (IVFS for short)  $A$,  namely  a map
$A: Q \to L([0,1)]$ that assigns to each position  $(i,j)$ an interval $ \xb^{i,j}$ (see next formula \eqref{intorno}). \\

Let  $\alpha \in [0,1]$ be fixed, and let 
  $K_\alpha : L[0,1] \rightarrow [0,1]$  
be a function,  given in \cite{atanassov,Burillo1,Burillo2,couso}, such that 
for every  $\xb \in L([0,1])$ and $\alpha \in [0,1]$,
\begin{description}
\item[k.1)] 
$ K_0 ( \xb)=\xg $, \quad $ K_1( \xb )=\xs,$  \quad $K_{\alpha} (\xb) = \xg$ if  $\xg=\xs$;
\item[k.2)] for every  $ \alpha \in [0,1] \,\,  K_\alpha(\xb)= K_0(\xb) + \alpha(K_1(\xb) -K_0(\xb))$;
\item[k.3)] if $\xb \leq_L \yb,\, \xb, \yb \in L([0,1])$ then  $K_{\alpha} (\xb) \leq K_{\alpha} (\yb)$ for every $\alpha \in [0,1]$;
\item[k.4)] $ \alpha \leq \beta$ if and only if  $K_{\alpha} (\xb) \leq K_{\beta} (\xb)$ for every $\xb \in L([0,1]).$
\end{description} 
The operator $K_\alpha$ is known in the literature as  Atanassov's operator.\\

Using $K_\alpha$
 it is  possible to associate  an interval-valued fuzzy set  with a fuzzy set in  the following way:
\begin{eqnarray}\label{opA}
K_\alpha(\xb)= K_\alpha( [\xg,\xs])=\xg + \alpha(\xs-\xg) = \xg + \alpha W(\xb).
\end{eqnarray}
In practice,  Atanassov's operator of order $\alpha$ is a convex combination of the end points of its argument
 $\xb=[\xg, \xs]  \in L[0,1]$.

\begin{remark}\label{rem1} \rm
There are other possible constructions of the multifunction $K_{\alpha}$,  and
 the choice of the previous operator  is motivated by the length of the interval 
 being  fundamental in the magnification process given in \cite{Mesiar}, since the length of each interval membership  is  fixed a priori.	
\end{remark}	

\subsection{Interval-valued fuzzy model}\label{IVFM}

We  provide  a  description of the algorithm based on the above interval-valued fuzzy model.
 For the sake of brevity,  we often refer to such an algorithm  with the term "fuzzy-type algorithm".\\

As previously mentioned  let $Q$  be an  $n \times m$ matrix associated  with a grayscale image.
Let $\delta \in [0,1]$ and $p \in \mathbb{N}$.
For every  $i \in \{1,2, \ldots,m\}$ and $j \in \{1,2, \ldots, n\}$ let  $q_{i,j}$ be the value of \sout{the} element $(i,j)$  in  $Q$.\\
For every $i \in \{1,2, \ldots,m\}$ and  $j \in \{1,2, \ldots, n\}$
let $V_{i,j}=(v^{(i,j)}_{k,l})_{k,l}$ a $(2p+1) \times (2p+1)$  square matrix (also named block) centered at the position $(i,j)$,
 namely the value $v^{(i,j)}_{p+1,p+1}$ coincides with $q_{i,j}$, and  is used  to obtain the magnification of $Q$. \\
Let  $v^{(i,j)}_{k,l}$  be the elements of  $V_{i,j}$ with  $k,l \in \{1,2, \ldots, 2p+1\}$;
\begin{eqnarray}\label{vkl}
 v^{(i,j)}_{k,l} = \left\{ \begin{array}{ll}
q_{i-p+k-1,j-p+l-1} & \mbox{if } i-p+k-1 \in \{1,2, \ldots, n\} , \\
& \phantom{aa}  j-p+l-1 \in \{1,2, \ldots, m\} \\
0 & \mbox{elsewhere.}
\end{array} \right. 
\end{eqnarray}
This means that if  there are positions $(k,l)$  in $V_{i,j}$ that are not covered by elements of $Q$
 (i.e., if in  the superposition  of the block $V_{i,j}$
 with  the matrix $Q$ there are some elements that do not belong to $Q$), the corresponding values  in the matrix $V_{i,j}$ 
are  set to zero.
\\
 To define a   neighborhood of  $q_{i,j}$,  the oscillation $\omega_{i,j}$  of the values in $V_{i,j}$
   is calculated  without considering the possible presence of the added null  values in the block , namely,
\begin{eqnarray}\label{oij}
\omega_{i,j} &=& \left(  \max_{\substack{i-p+k-1 \in \{1,2, \ldots, n\},\\ j-p+l-1 \in \{1,2, \ldots, m\}}}
 q_{i-p+k-1,j-p+l-1} \right) + \nonumber \\
&-& \left( \min_{\substack{i-p+k-1 \in \{1,2, \ldots, n\},\\ j-p+l-1 \in \{1,2, \ldots, m\}}}
 q_{i-p+k-1,j-p+l-1}\right),
\end{eqnarray}
 and a closed interval
  $F(q_{i,j},\omega_{i,j}, \delta) \in L([0,1])$ is assigned to each $q_{i,j}$,  as follows:
\begin{eqnarray}\label{intorno} 
F(q_{i,j},\omega_{i,j}, \delta) = [ q_{i,j} ( 1 - \delta \omega_{i,j}), q_{i,j} ( 1 - \delta \omega_{i,j}) + \delta \omega_{i,j}].
\end{eqnarray}
 Therefore the intensities of the pixels in this generated block provide  information for  
 obtaining the length of the interval-valued membership  built using $F$. 
 For  this interval-valued membership   in $L([0,1])$,   Atanassov's operator
\eqref{opA}  is applied  to construct a new square matrix
\[V'_{i,j}=(v'_{k,l})_{k,l},  \qquad k,l \in \{1, 2, \ldots, 2p+1\},\] 
whose elements are obtained in the following way: 
\begin{eqnarray*}
 v'_{k,l} &:=& K_{v^{(i,j)}_{k,l}}(F(q_{i,j}, \omega_{i,j}, \delta ) )  
=K_{v^{(i,j)}_{k,l}} ([ q_{i,j} ( 1 - \delta \omega_{i,j}), q_{i,j} ( 1 - \delta \omega_{i,j}) + \delta \omega_{i,j}])\\
&=& 
v^{(i,j)}_{k,l} \cdot \big( q_{i,j} ( 1 - \delta \omega_{i,j}) + \delta \omega_{i,j} \big) 
+ (1 - v^{(i,j)}_{k,l} ) \cdot 
q_{i,j} ( 1 - \delta \omega_{i,j}) = \\
&=& 
v^{(i,j)}_{k,l}  \delta \omega_{i,j} + q_{i,j} ( 1 - \delta \omega_{i,j}).
\end{eqnarray*}

Finally, in the new rescaled image, each element 
 $q_{i,j}$ is replaced by the new block $V'_{i,j}$.
We  can observe that if $\delta = 0$ 
the information on the boundary is lost since 
$F(q_{i,j},\omega_{i,j}, \delta) = q_{i,j}$.
\\

\section{Other methods}
 To evaluate the performance of the considered fuzzy-type algorithm, in the numerical tests performed in 
Section \ref{sec5},  we  consider the rescaling of a given dataset of images with the well-known bicubic method, 
 which is very classical in digital image processing, and  is already implemented in several software and 
 dedicated commands are available in most used programming languages) and we compare it  with the SK algorithm 
which will be recalled in the next subsection. 

\subsection{The Sampling Kantorovich algorithm for image rescaling}
An algorithm  that has been widely applied in the field of image rescaling is 
 known  for its name, the sampling Kantorovich (SK) algorithm; see, e.g., \cite{bbsv,Poland1}.
 The above tool arises as an optimized implementation of a family of sampling-type operators, that is, the 
 multivariate SK operators, defined through the following formula:
\begin{equation} \label{KANTOROVICH}
(S_w f)(\vec{x})\ :=\ \sum_{\vec{k} \in \mathbb{Z}^{2}} \chi(w\vec{x}-\vec{k})\, 
\left[ w^2 \int_{R_{\vec{k}}^w }f(\vec{u})\ d\vec{u} \right], \hskip0.7cm \vec{x} \in \mathbb{R}^2, \quad w>0,  
\end{equation}
where $f: \mathbb{R}^2 \to \mathbb{R}$ is a locally integrable function (signal/image)
 such that the above series is convergent for every $\vec{x} \in \mathbb{R}^2$, and 
\[
R_{\vec{k}}^w\ :=\ \left[\frac{k_1}{w},\frac{k_1+1}{w}\right]\times\left[\frac{k_2}{w},\frac{k_2+1}{w}\right],
\]
are the squares in which we consider the averaged values of the sampled signal $f$ (see for example \cite{CANCV4,CACOVI6}).

 $S_w$, $w>0$, are approximation operators 
 that can   pointwise reconstruct continuous and bounded signals, and to uniformly reconstruct
 signals  that are uniformly continuous and boun\-ded, as $w \to +\infty$. 
Moreover, the $S_w$ operator  can also be used  to reconstruct not-necessarily continuous signals, e.g., signals belonging
 to the $L^p$-spaces, $1 \leq p < +\infty$ (\cite{AT,AV,COVI,cma,MIMMO,carlo,bbsv,cand,carlo0,BS1,BS2,KS1,KS2,VZ,CS2018,CSV2019}).
  The function $\chi: \mathbb{R}^2 \to \mathbb{R}$,  given in \eqref{KANTOROVICH}, 
 is called a {\it kernel} and it satisfies  the following suitable assumptions,  very typical in this situation,  which
 are the usual conditions assumed by  discrete approximate identities (for more details, see, e.g., \cite{ACDR1}).
 Below, we present a list of functions that can be used as kernels in the formula recalled in \eqref{KANTOROVICH}.
      
      First, we recall the definition of the one-dimensional central B-spline of order $N$ (for example see  \cite{BSS1}):
\begin{equation} \label{splines}
 \beta^N(x)\ :=\ \frac{1}{(N-1)!} \sum^N_{i=0}(-1)^i \binom{N}{i} 
       \left(\frac{N}{2} + x - i \right)^{N-1}_+,     \hskip0.5cm   x \in \mathbb{R}.
\end{equation}
The corresponding bivariate version of  the  central B-spline of order $N$ is given by:
\begin{equation}
{\mathcal B}^N_2(\vec{x})\ :=\ \prod^2_{i=1}\beta^N(x_i), \hskip1cm \vec{x}=(x_1,x_2) \in \mathbb{R}^2.
\end{equation}

  Other important kernels are given by the so-called Jackson type kernels of order $N$, defined in the univariate case by:
\begin{equation}\label{jack}
J_N(x)\ :=\ c_N\, \mbox{sinc}^{2 N}\left(\frac{x}{2 N\pi}\right), \hskip1cm x \in  \mathbb{R},
\end{equation}
with $N \in \mathbb{N}$ and $c_N$ is a nonzero normalization coefficient, given by:
\[
c_N\ :=\ \left[ \int_{\mathbb{R}} \mbox{sinc}^{2 N}\left(\frac{u}{2 N \pi } \right) \, du \right]^{-1}.
\]
For the sake of completeness, we recall that the well-known $sinc$-function is 
 defined as $\sin (\pi x)/\pi x$, if $x \neq 0$, and $1$ if $x=0$; see e.g., \cite{KS1,KS2}. 
As in the  case of the central B-splines,  the bivariate Jackson type kernels of order $N$ are defined by:
\begin{equation}
{\mathcal J}^2_N(\vec{x})\ :=\ \prod^2_{i=1}J_N(x_i), \hskip1cm \vec{x}=(x_1,x_2) \in \mathbb{R}^2.
\end{equation}
In particular, Jackson type kernels have been
 revealed to be very useful, e.g., for applications  in   the biomedical field, \cite{Poland1}. 
For the numerical tests given in this paper, we  consider the bivariate Jackson-type kernel 
 with $N$ varying from 2 to 12. This choice will be motivated later.
 For several examples of kernels, see, e.g., \cite{NAT1,CPV1,CPV3,CNV2};  for more details about the SK
 operators and the corresponding SK algorithm, see e.g., \cite{CLCOGUVI1}, where  a pseudo-code
 is also available. For some applications of the SK algorithm to real world problems involving images, see, e.g., \cite{CSV3,Poland1}.

\section{Comparisons and evaluation of the numerical results: likelihood index $S$ and PSRN}

 To compare  the considered algorithms for image rescaling, we 
 use the following  indices that are known in the literature. 

The first  tool is  the Peak Signal-to-Noise Ratio (PSNR),  which is a well known 
index in the literature and  is often used to quantify the rate of similarity between two general signals. 

 The PSNR is defined  as  the Mean Square Error (MSE):
\[MSE\ =\ \sum\limits_{i=1}^N \sum\limits_{j=1}^M \dfrac{\left|I(i,j)-I_r(i,j)\right|^2}{NM},\]
where $I$ is the original image, $I_r$ is the reconstructed version of the original image $I$, $N$ and $M$ are the dimensions of the images. 
 Therefore  the PSNR is generally defined as follows:
\[
PSNR\ =\ 10\cdot \log_{10}\left(\dfrac{f_{max}^2}{MSE}\right),
\]
where $f_{max}$ represents the maximum value of the considered pixel's scale. 
 For  $8$-bit gray scale images  $f_{max}=255$, while  for  images 
 with  pixel values between $0$ and $1$ ( such  as those considered in our fuzzy algorithm)  $f_{max}=1$. 

     Hence, the PSNR formula  used in this paper  is expressed as follows:
\begin{equation}\label{psnr}
PSNR\ =\ 10\cdot \log_{10}\left(\frac{1}{MSE}\right).
\end{equation}

It is clear from the above definition that, the similarity between two images 
 is greater for the highest values of the PSNR.
  
   Furthermore, we  use another useful similarity index,   called the likelihood index S,  
 which was introduced by Bustince, et al. (\cite{BPB2007}), and  is defined as follows:
\begin{equation} \label{index-S}
S\ :=\ {1 \over N \times M}\, \sum_{i=1}^N \sum_{j=1}^M \left[ 1 - \left|I(i,j)-I_r(i,j)\right|  \right],
\end{equation}
where the notations used in (\ref{index-S}) are the same as those  employed in the definition of the PSNR \eqref{psnr}. 
It is clear from the above definition that, the parameter $S$ can assume values between $0$ and $1$, and that for
 closer images $S$ should be as close as possible to $1$. 

\section{Numerical experiments} \label{sec5}

In this section we provide a numerical comparison among the algorithms considered in the previous sections, namely
 the fuzzy-type algorithm, the classical bicubic  and   the SK  algortithm. 
Such a comparison will be carried out thanks to the similarity 
 indices previously recalled, i.e., the PSNR and the likelihood index S.

For the numerical tests, we proceed as follows. 
  We first consider a set of original images of a given dimension $N \times M$, 
 which  will be used as a reference. 
Such images will be reduced without interpolation (using the nearest neighbor method \cite{BIAU}) to the dimension
 ${N \over 3} \times {M \over 3}$.\\
 Finally, the reduced images will be rescaled to the original dimension by using the methods mentioned above.
In this way we  dispose of a reference image (the original image), and three reconstructed images generated by
 the three different methods  mentioned above. \\
With respect to
 the application of the algorithm based on sampling Kantorovich operators, in view of the accurate experimental analysis
 given in \cite{CSV3},
 the SK algorithm has been applied using the parameters that have been seen to be the best possible under certain
 qualitative criteria (for more details see \cite{CSV3} again). 
More precisely, we  consider  the
 bivariate Jackson-type kernel 
 ${\mathcal J}^2_{N}$ with $N \in \{2, 3, \ldots 12\}$.\\
 Concerning the parameter $w$ in the SK algorithm, we consider the following values:
 $w=5, 10, 15, 20$,
 and 
  $25$ only for the baboon image.\footnote{Note that, as stated in \cite{CSV3}, in the case of the rescaling of images with double dimensions,
 it is sufficient to choose $w=15$ when $N =12$.} \\
 
The image dataset (the source files are contained in the repository  
 https://links.uwa\-ter\-loo.ca/\-Repository.html or in \cite{foto})  is composed   of the four different 
 grayscale images shown in Figure \ref{fig1}. There are  the classical "baboon" and "boat", 
 which are  commonly used in image analysis, and two pictures of a "city" and a "mountain", respectively.\\
{
\begin{figure}[h!]
\hskip-.1cm \includegraphics[scale=0.52]{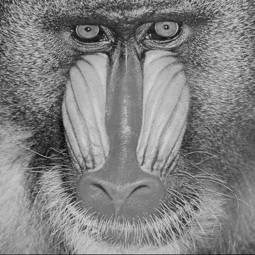}
\hskip0.3cm
\includegraphics[scale=1.1]{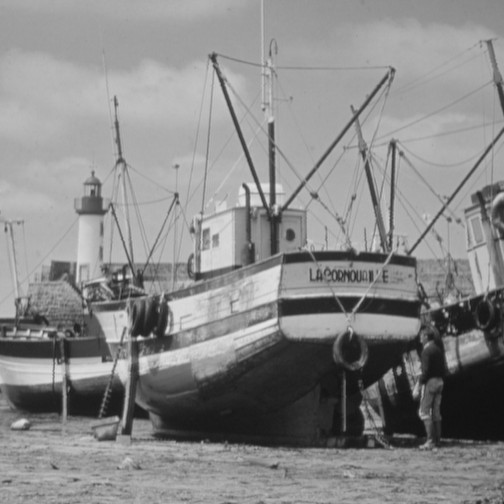}
\vskip0.4cm
\includegraphics[scale=0.66]{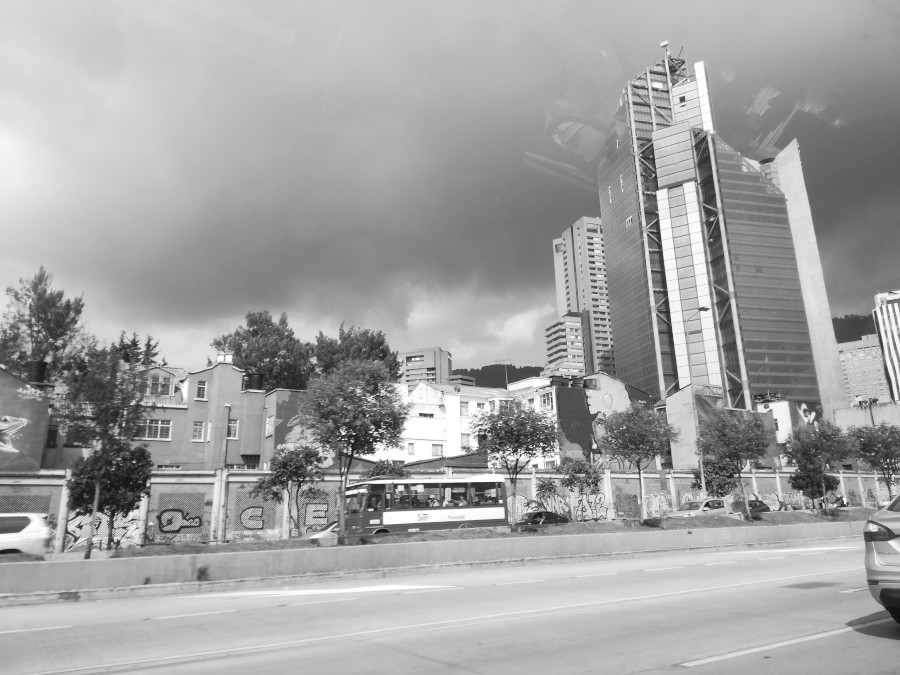}
\hskip0.35cm
\includegraphics[scale=1]{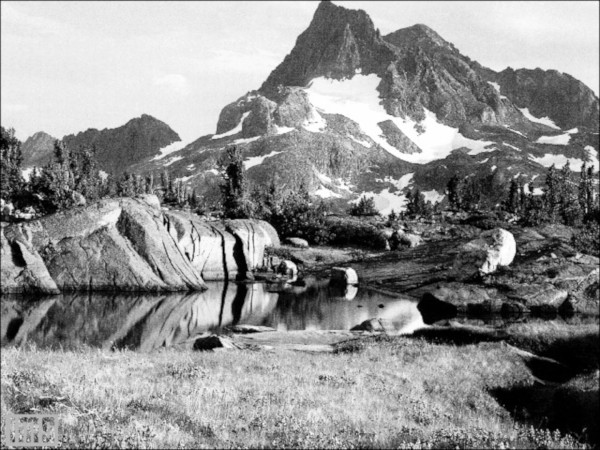}
\caption{\footnotesize Reference images: baboon ($255 \times 255$ pixel resolution); boat ($504 \times 504$ pixel resolution); 
city ($675 \times 900$ pixel resolution); mountain ($450 \times 600$ pixel resolution).} \label{fig1}
\end{figure}
}
\\
 The choice of the four images  is motivated 
by the fact that we want to  compare  images
of different sizes, brightness levels  and textures.
Finally  the boat image was also considered in the quoted paper \cite{Mesiar}, but we do not know if it has the same
 dimension or resolution.
The histograms of the four images
show  that the distributions of the grayscale of the various images are very  different from each other.
{
\begin{figure}[h!]
\hskip-.1cm \includegraphics[scale=0.33]{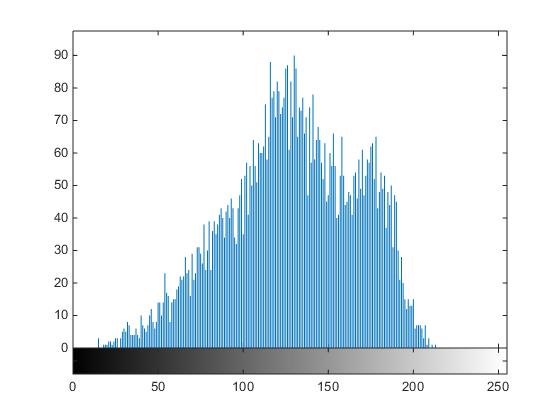}
\hskip0.2cm
\includegraphics[scale=.33]{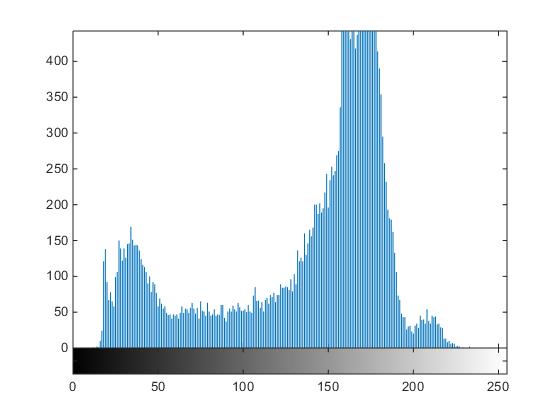}
\vskip0.5cm
\includegraphics[scale=0.33]{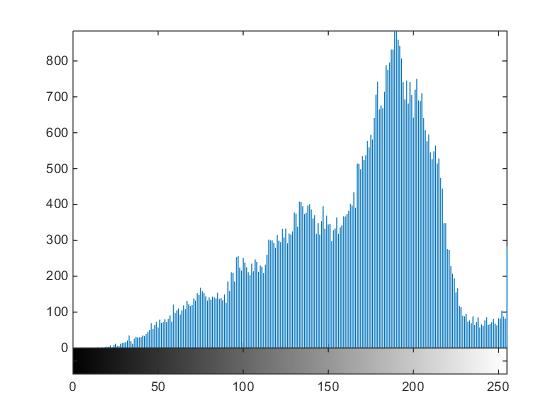}
\hskip0.2cm
\includegraphics[scale=0.33]{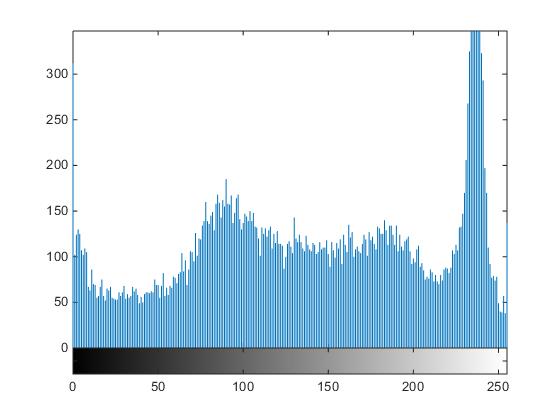}
\caption{\footnotesize Histograms of the original images in the dataset: baboon, boat,  city, mountain.} \label{fig-hist}
\end{figure}
}

 The empirical simulation of the two algorithms is performed on Windows 11 operating system
   with an Intel Core i7 8th gen. Moreover, all the programs
are written and compiled on MATLAB version R2014b.\\

Concerning the application of the fuzzy-type algorithm, we provide the rescaled images for values of the parameter
 $\delta$ running between $0$ and $1$,
 with a step-size equal to $0.01$, for each  of the images given in Figure \ref{fig1}. 
The corresponding results of the PSNR and likelihood index S  are plotted in Figure \ref{fig2}.
\newpage
%
{
\begin{figure}[h!]\label{delta-fig}
\fbox{\includegraphics[scale=0.262]{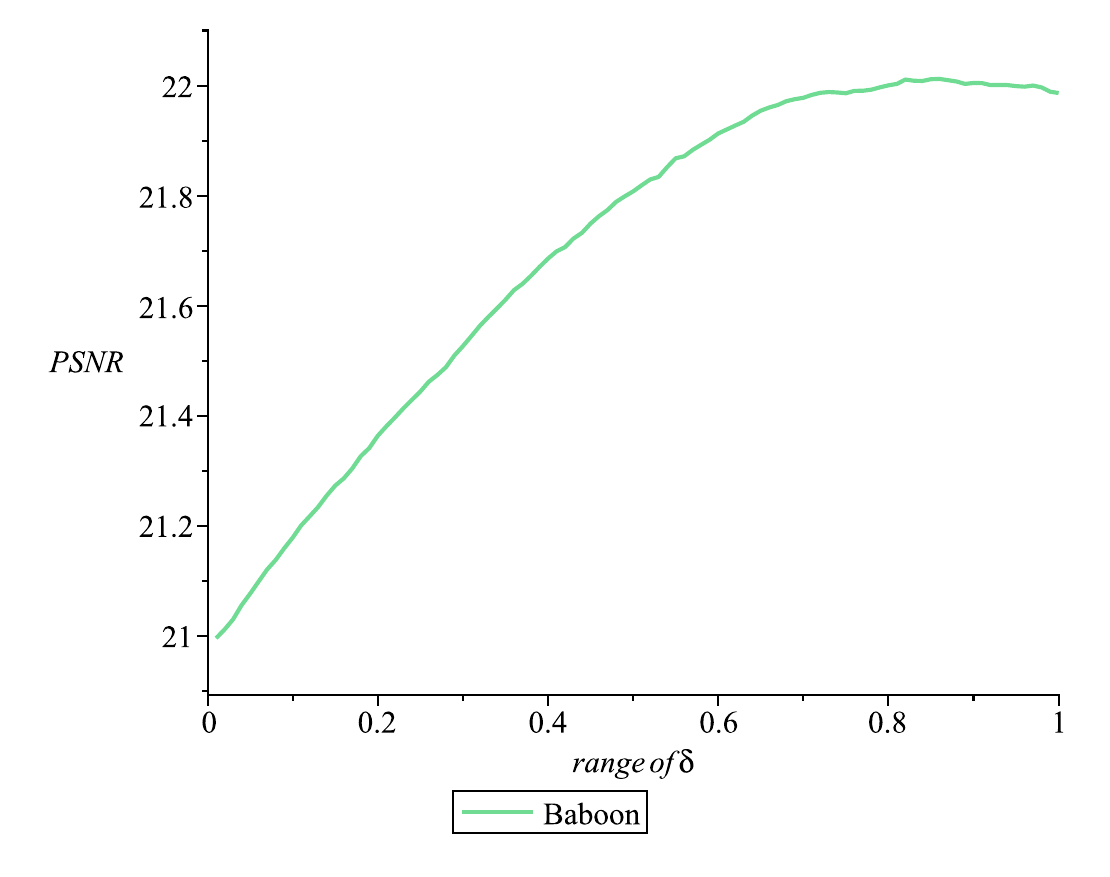}
\hskip0.3cm
\includegraphics[scale=0.262]{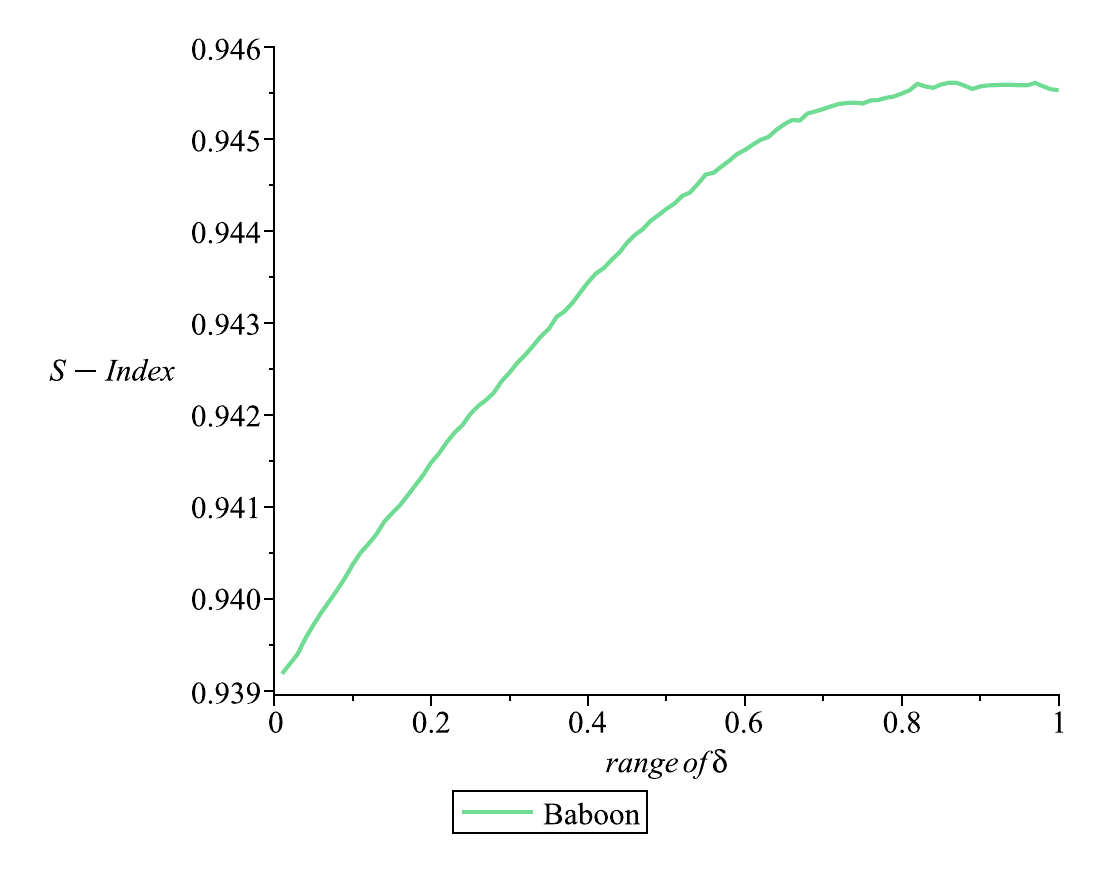} \quad }
\vskip0.4cm
\fbox{\includegraphics[scale=0.262]{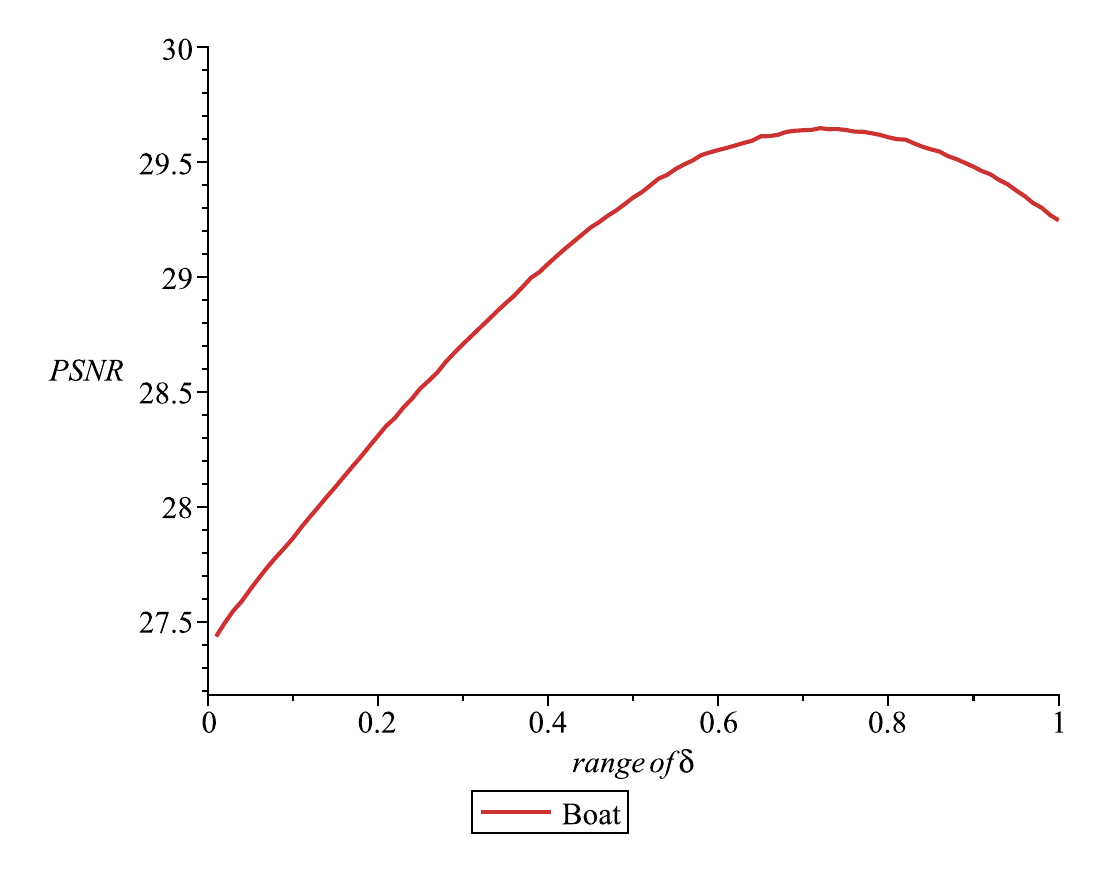}
\hskip0.3cm
\includegraphics[scale=0.26]{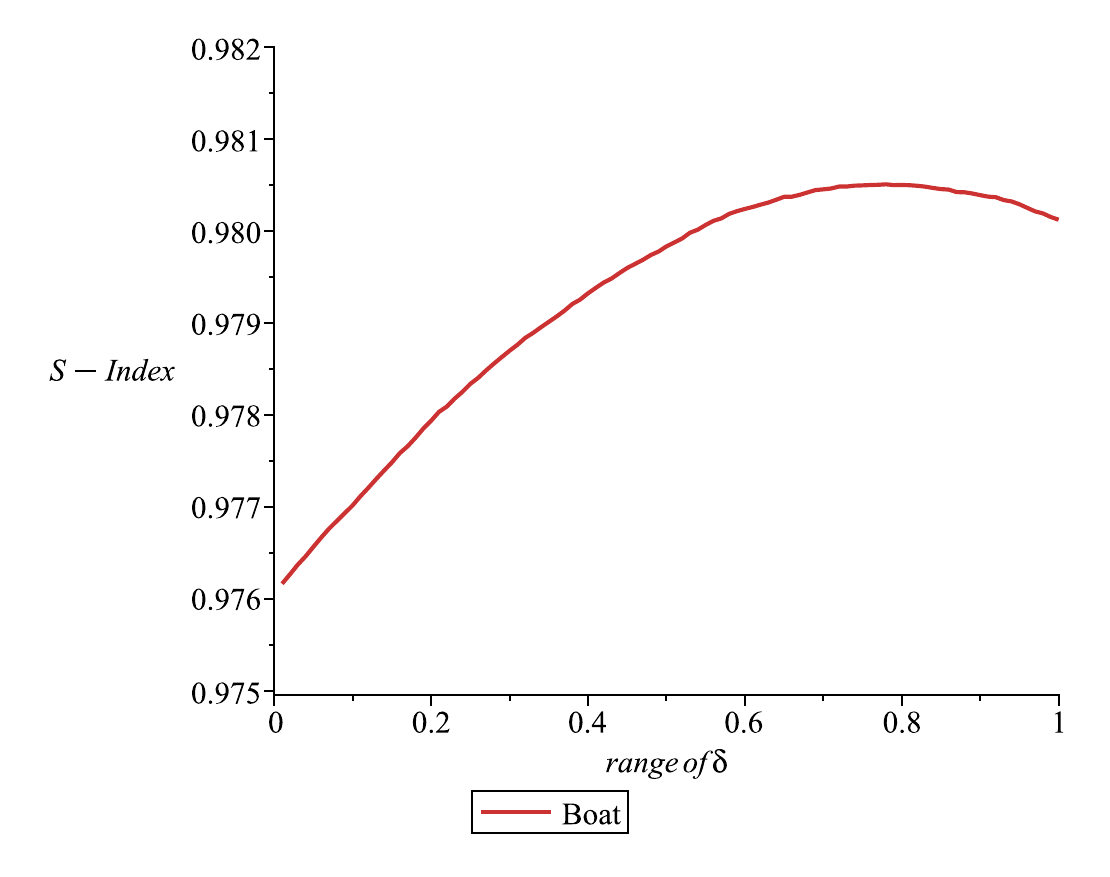}\quad }
\vskip0.4cm
\fbox{\includegraphics[scale=0.262]{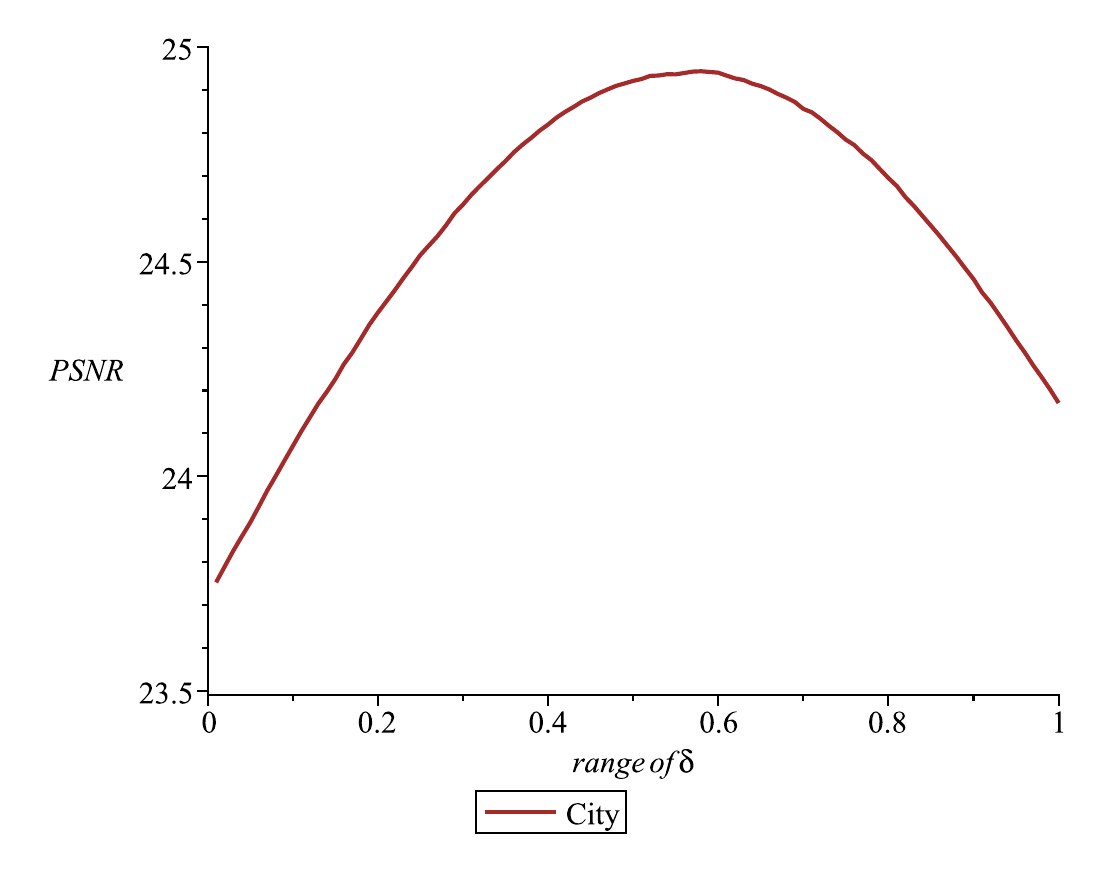}
\hskip0.3cm
\includegraphics[scale=0.262]{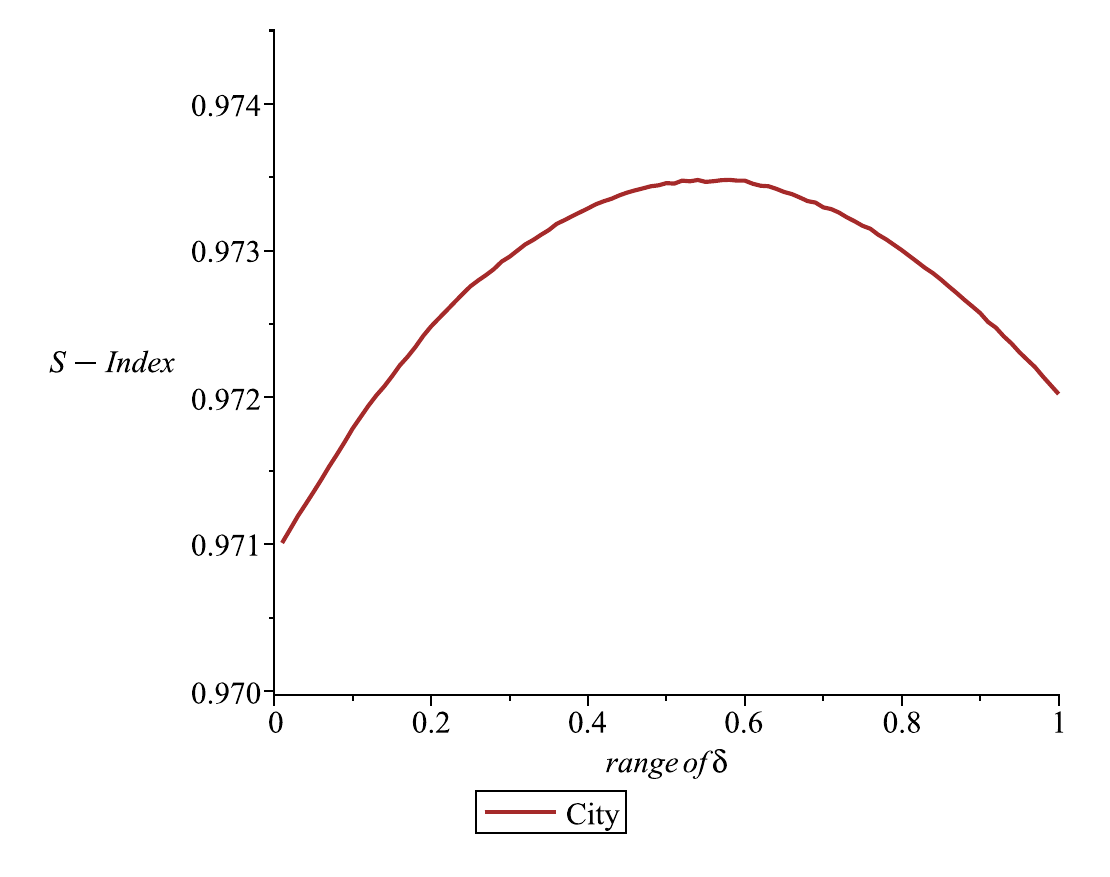} \quad}
\vskip0.4cm
\fbox{\includegraphics[scale=0.26]{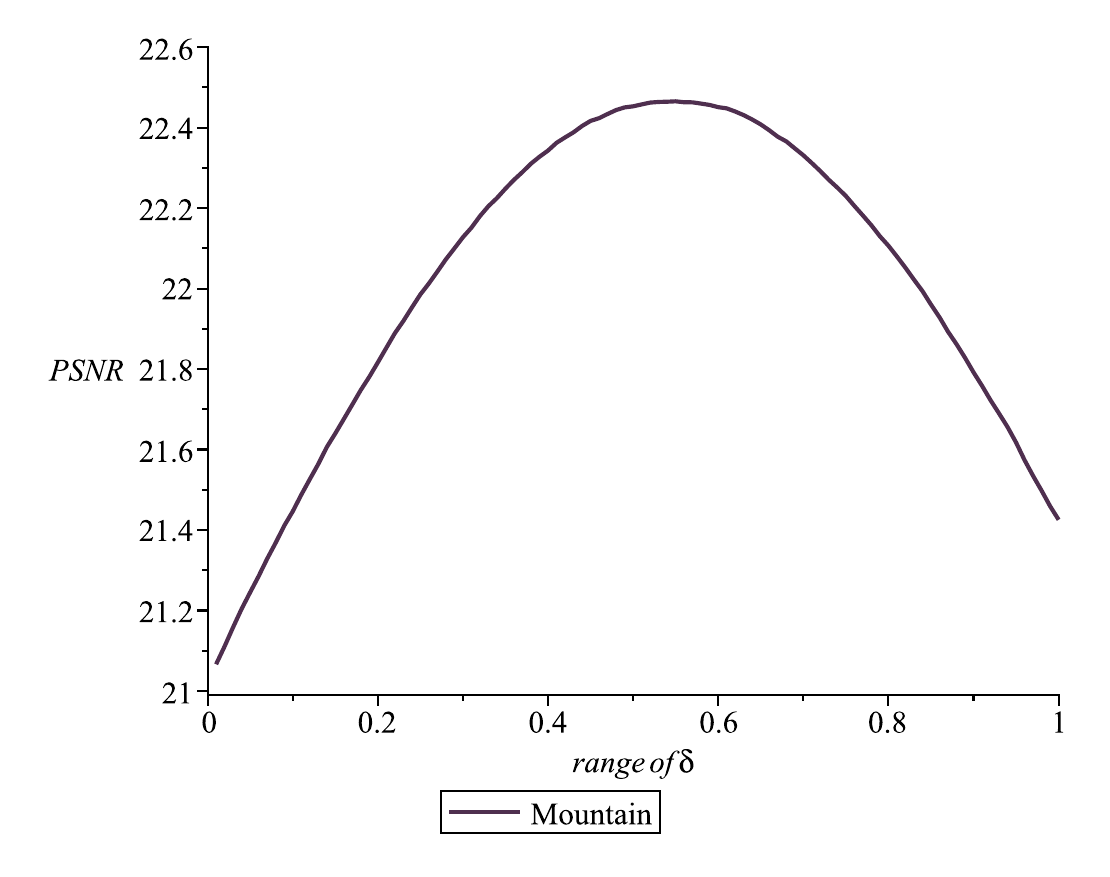}
\hskip0.3cm
\includegraphics[scale=0.26]{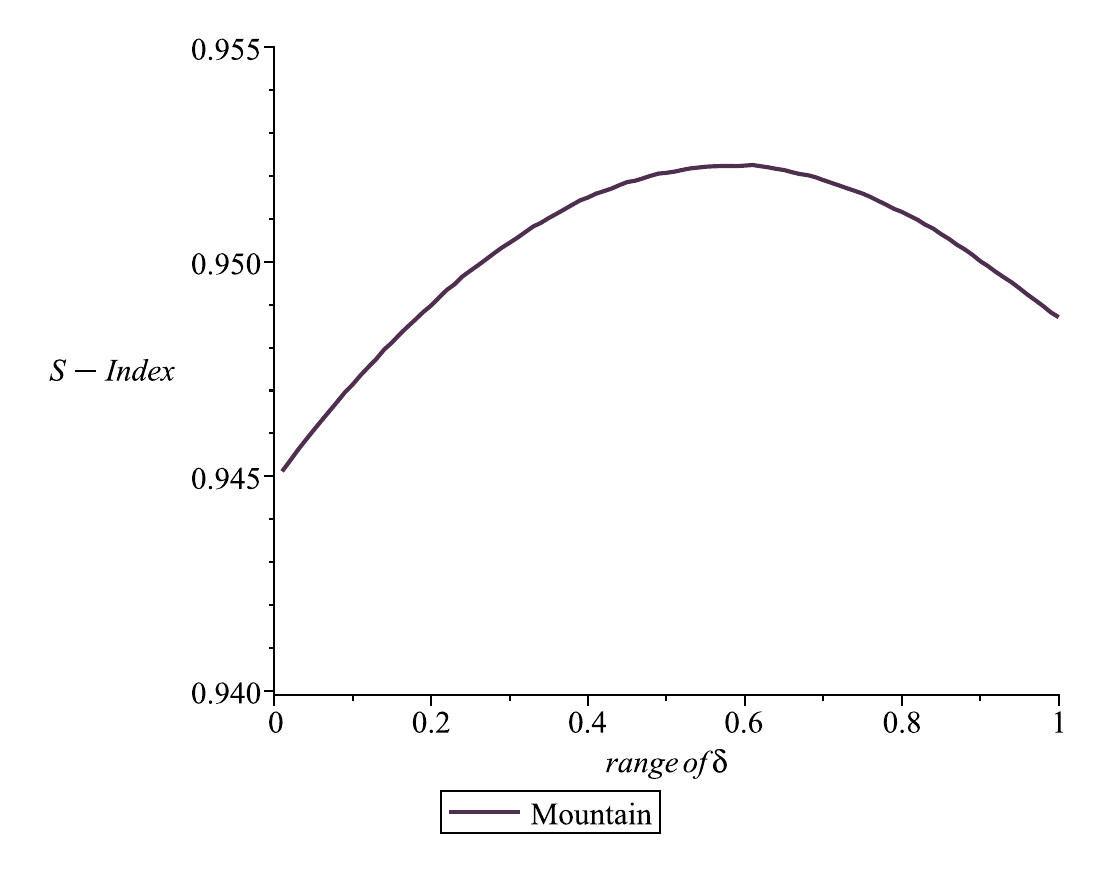} \quad}
\caption{\footnotesize The plots of the values of the PSNR and likelihood index S computed for the whole dataset of
 reconstructed images by the fuzzy-type 
algorithm when the parameter $\delta$  varies from $0$ to $1$ with step-size of $0.01$.} \label{fig2}
\end{figure}
}

Similarly we rescaled the dataset  images  using the Sampling Kantorovich algorithm and we  examined the values of 
the similarity indices corresponding to  the parameters $N$ and $w$ given in \eqref{KANTOROVICH} and \eqref{jack}. 
 The $w$ parameter determines the amount
 of the sample values that are involved in the reconstruction process, while $2N$ represents the order of decay of the considered kernel function.
 In particular, we examined  the parameter $N$ varying in the set $\{2,3, \ldots, 12\}$ and the parameter $w \in \{5,10,15,20,25\}$.
 Here there are
 plots of the values of the PSNR and likelihood  S indices computed for the  reconstructed dataset images with the Sampling Kantorovich
algorithm
for the considered values of 
 the parameter $N$ of the bivariate Jackson kernel  \eqref{jack}.
{
\begin{figure}[h!]
\fbox{\includegraphics[scale=0.25]{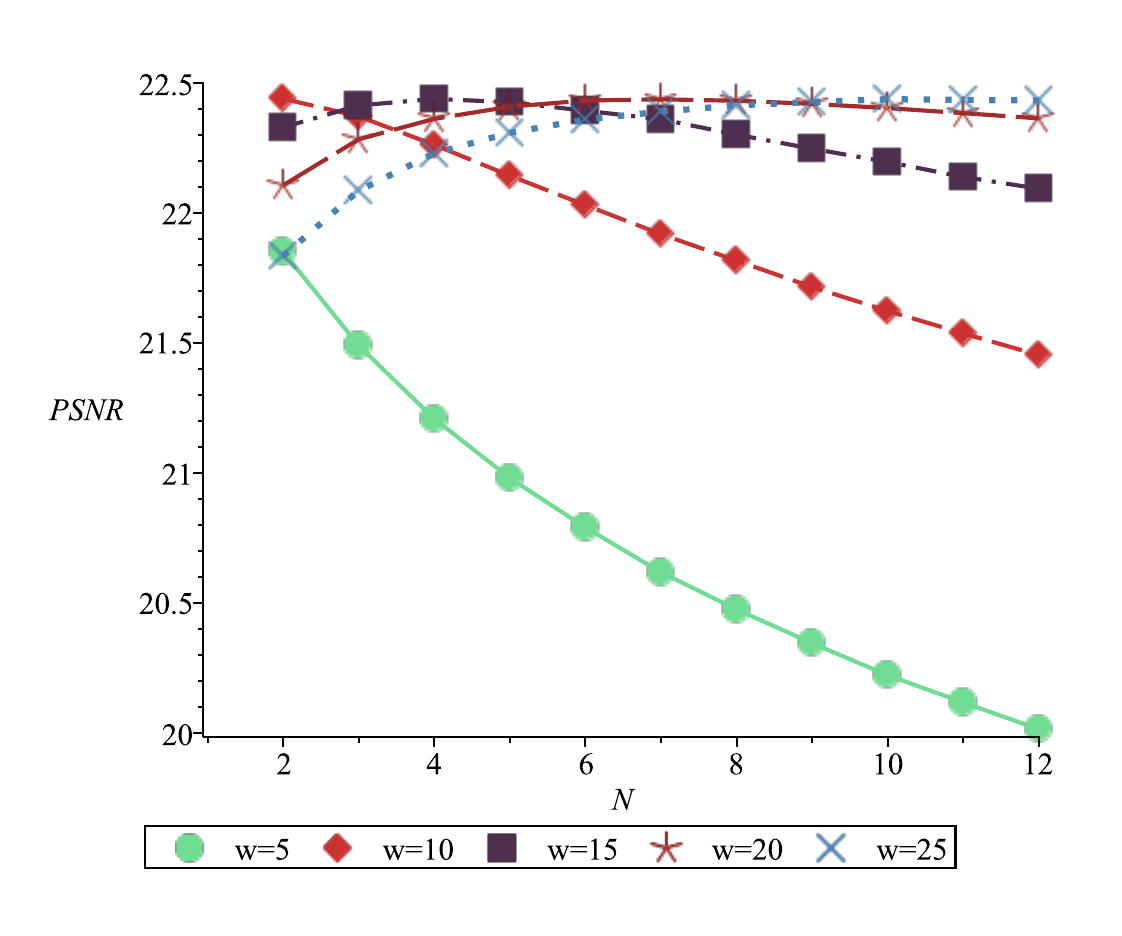}  \,\,  
\includegraphics[scale=0.25]{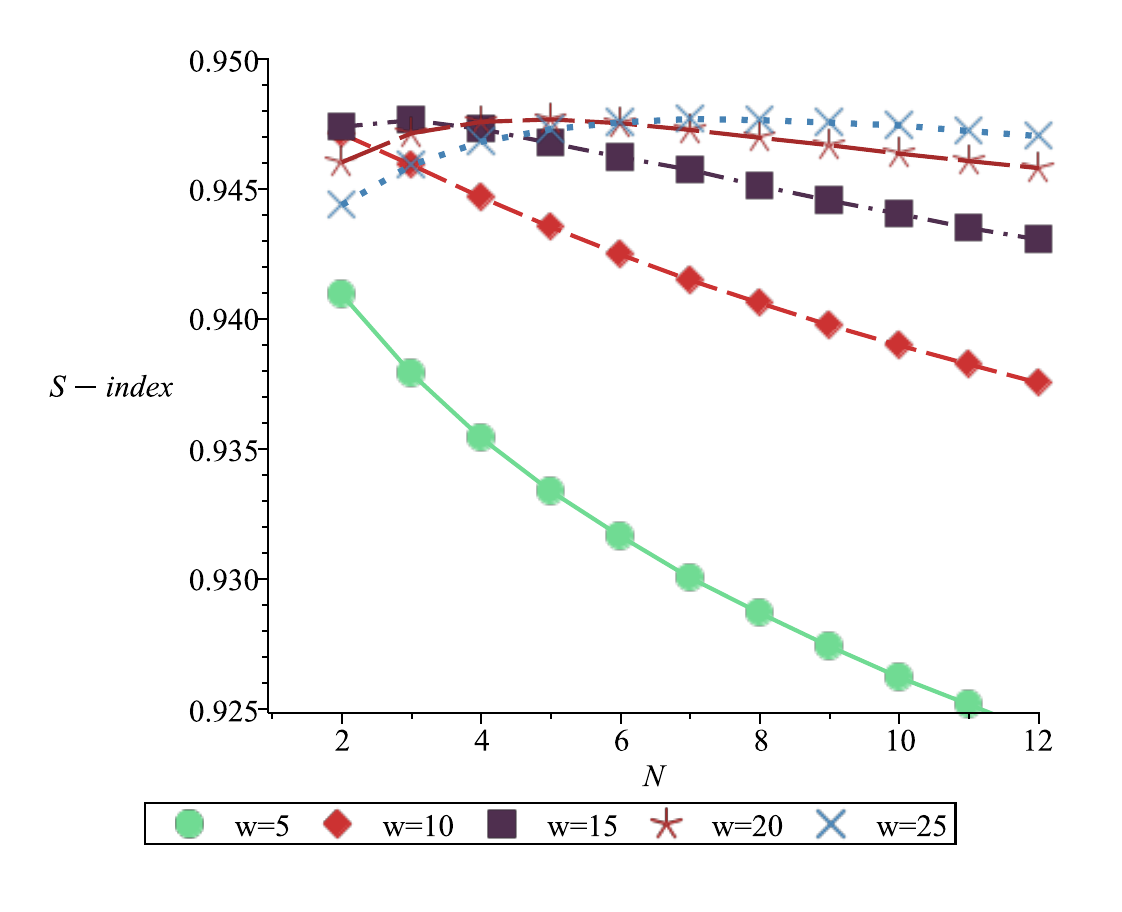} \quad}
\caption{\footnotesize The plots of the values of the computed indices for the  reconstructed baboon images.} \label{fig-baboon} 
\end{figure}
}

{
\begin{figure}[h!]
\fbox{\includegraphics[scale=0.25]{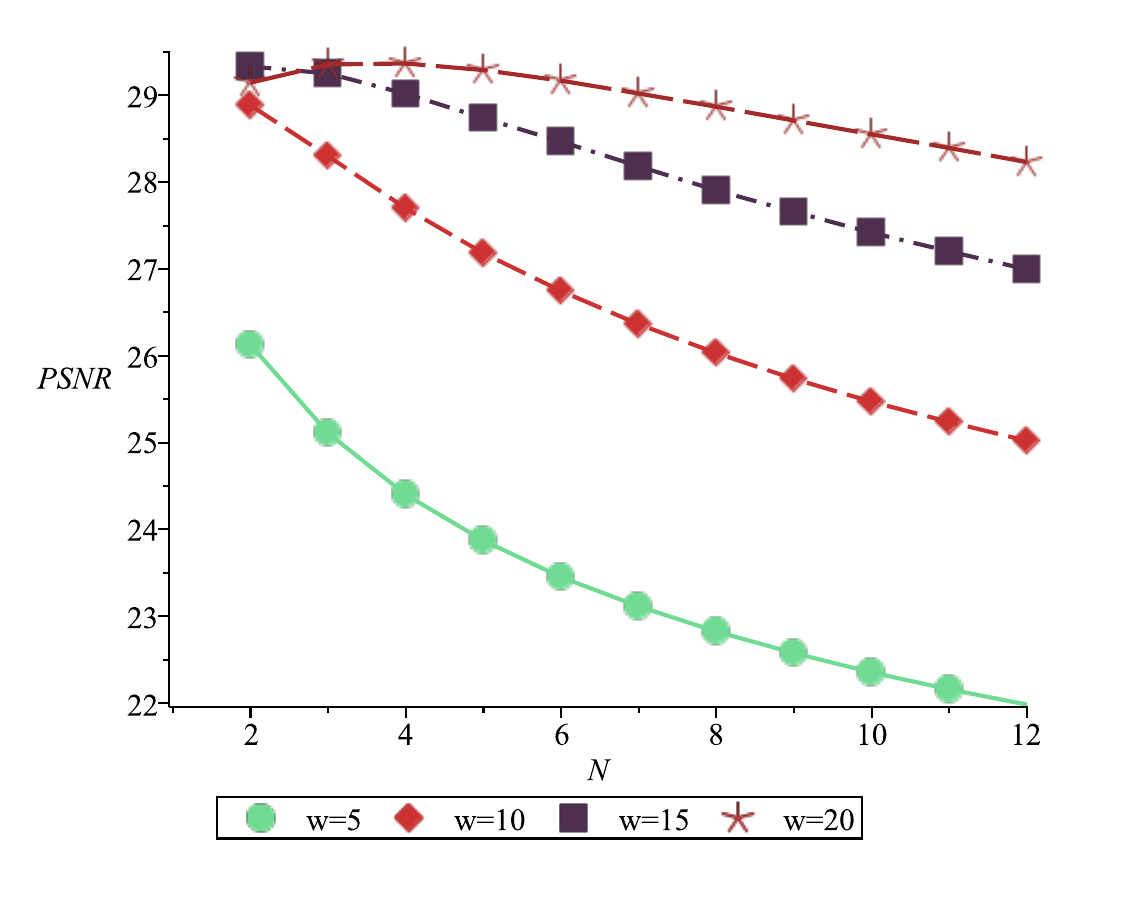}  \,\,  
\includegraphics[scale=0.25]{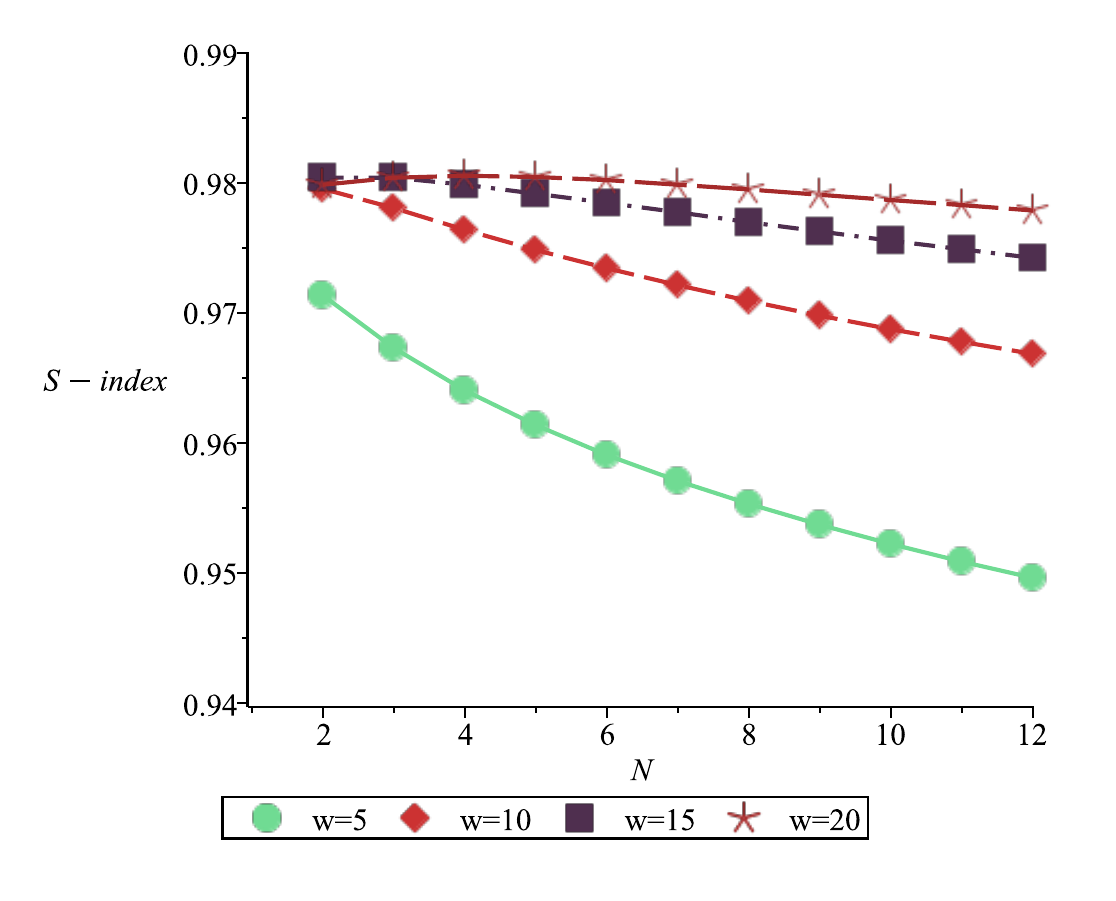} \quad}
\caption{\footnotesize The plots of the values of the computed indices for the  reconstructed boat images.} \label{fig-boat}
\end{figure}
}

{
\begin{figure}[h!]
\fbox{\includegraphics[scale=0.25]{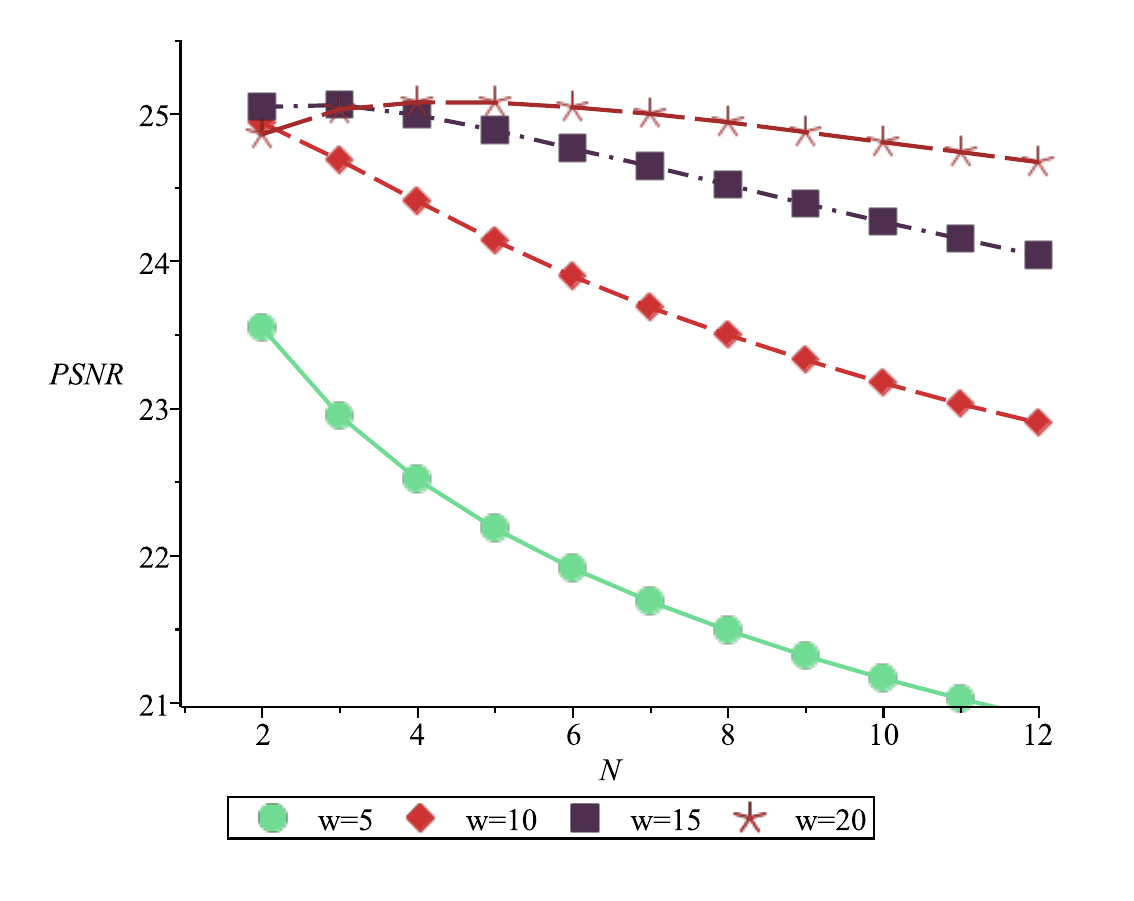} \,\, 
  \includegraphics[scale=0.25]{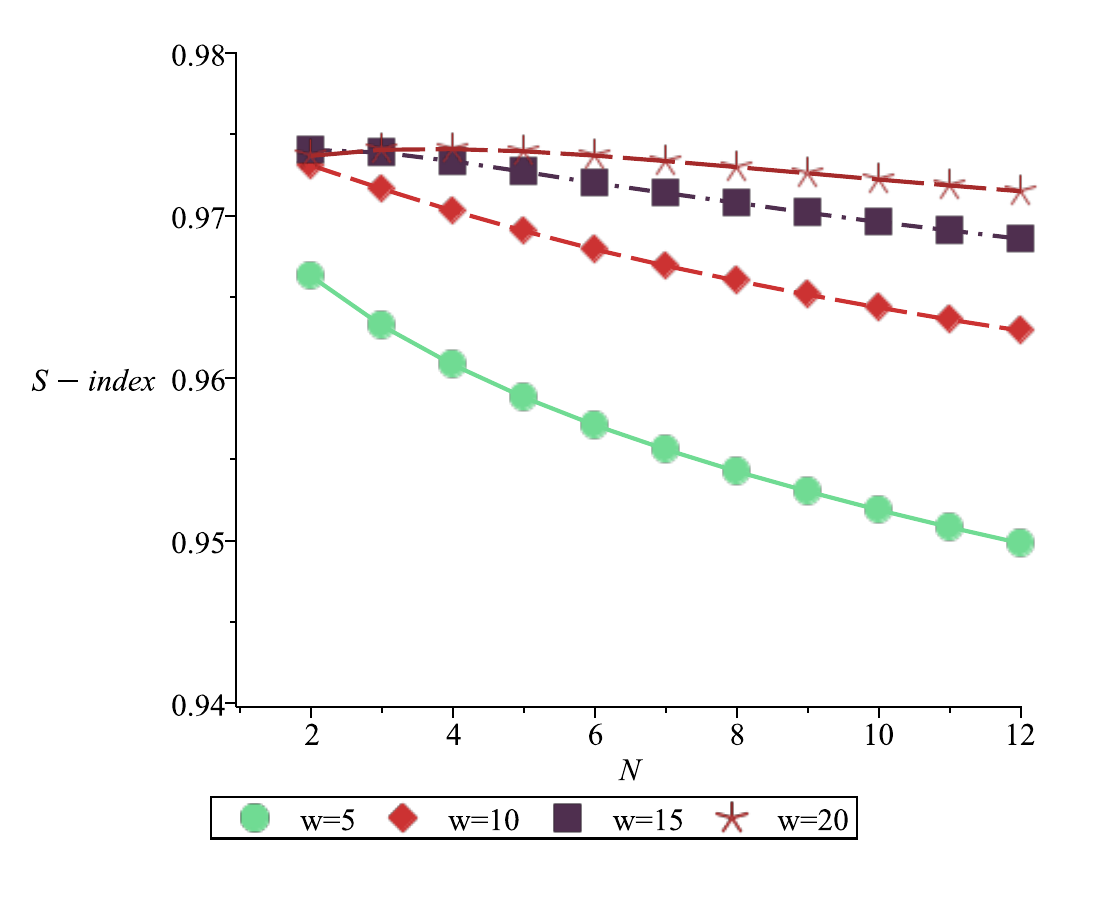} \quad}
\caption{\footnotesize The plots of the values of the computed  indices for the  reconstructed city images.} \label{fig-city}
\phantom{a} \vskip.24cm 
\end{figure}
}

{
\begin{figure}[h!]
\fbox{\includegraphics[scale=0.25]{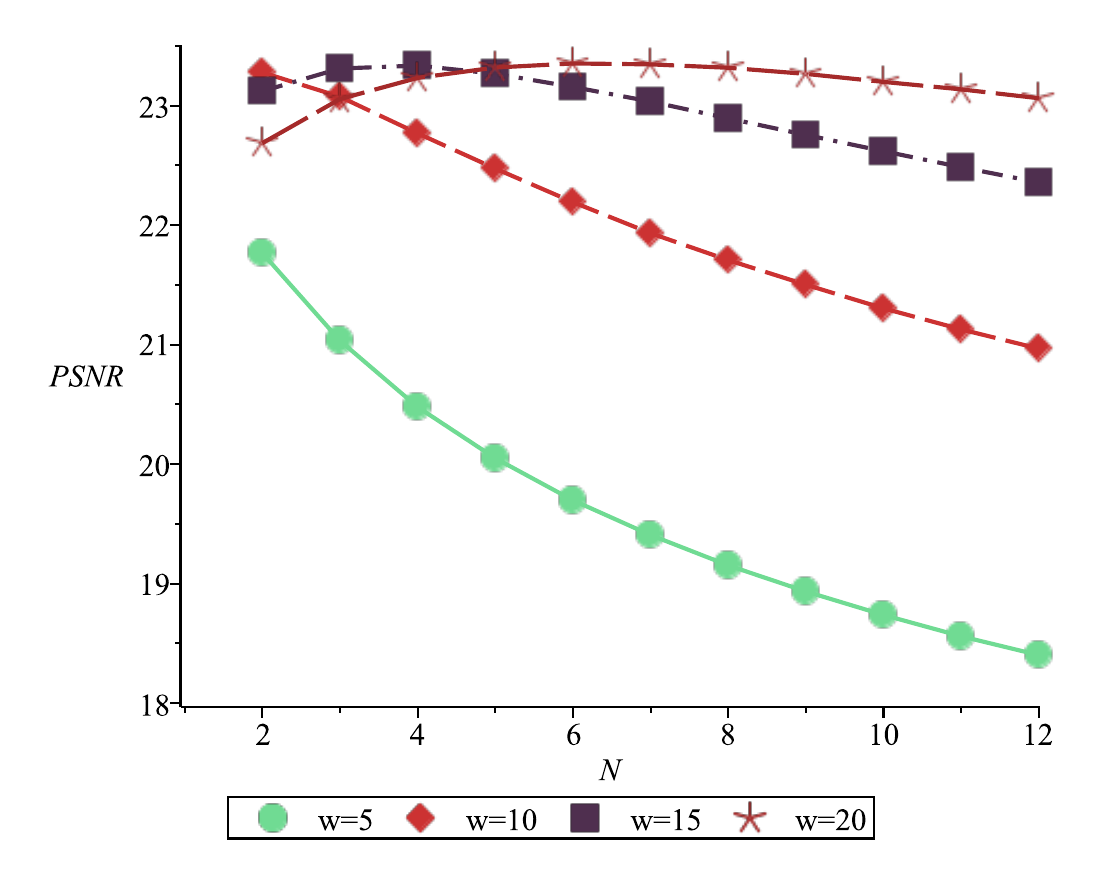}  \,\,  
\includegraphics[scale=0.25]{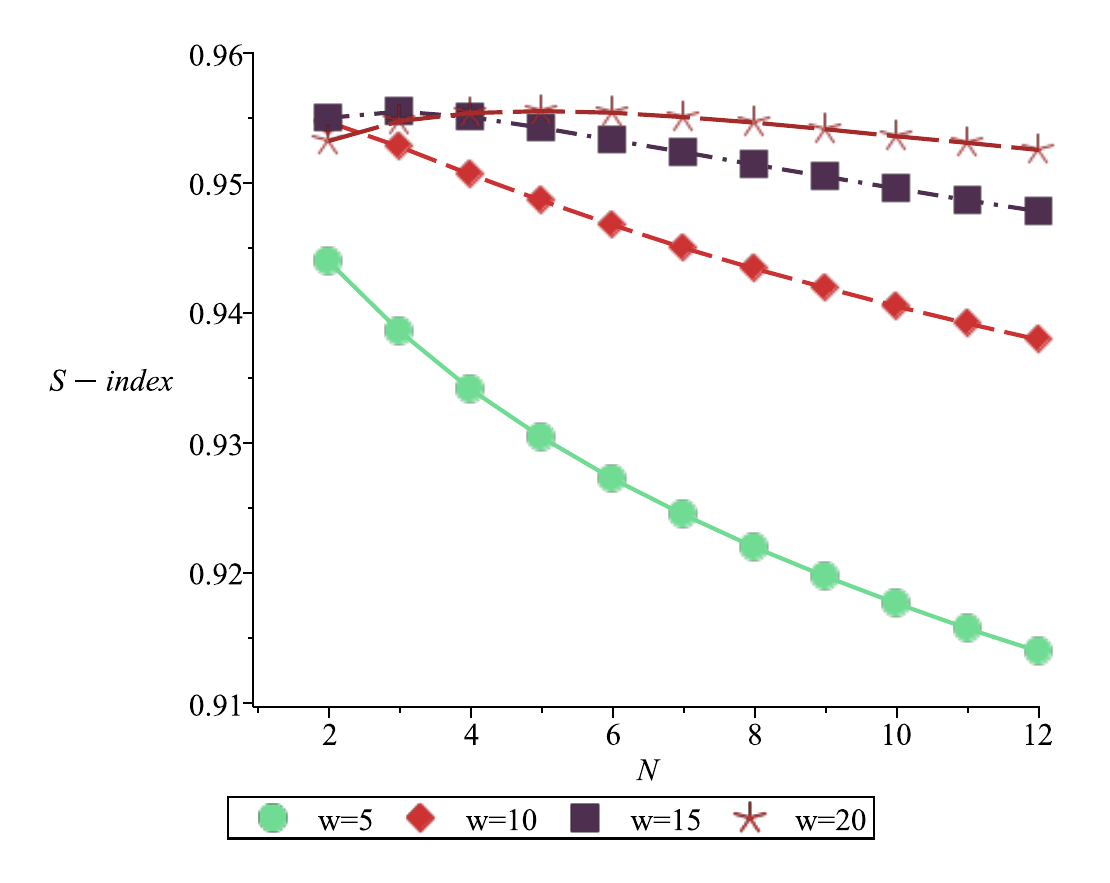} \quad}
\caption{\footnotesize The plots of the values of  the computed indices for the reconstructed mountain images.} \label{fig-mountain}
\end{figure}
}

For a more detailed version of these numerical outputs see  the Appendix  \ref{Appendix} where all indices are calculated for the SK algorithm.\\

To provide a more detailed comparison of the numerical results 
 shown in Figures \ref{fig2}, \ref{fig-baboon},
\ref{fig-boat}, \ref{fig-city} and \ref{fig-mountain} we also provide additional and useful data in the following tables.

  In Table \ref{tab1}, the values of the PSNR are listed and analysed in the case of the fuzzy-type algorithm,
 together with a comparison with the bicubic method and  the   SK algorithm.    In particular 
\begin{itemize}
\item  in the first column :
 "$\delta_{Psnr}^{\max}$"  denotes the  minimum value of the parameter $\delta$ for which the maximum PSNR is reached, 
when the images processed by the fuzzy-type algorithm are considered; 
\item "PSNR max - Fuzzy" denotes the maximum value of the PSNR reached  by  implementing
 the fuzzy-type algorithm for $ \delta^{\max}_{Psnr}$; 
\item  "PSNR - bicubic" denotes the values of the PSNR achieved by the image processed by the bicubic algorithm;
\item  ``$(N,w)_{Psnr}^{\max}$'' denotes the value of the pair $(N,w)$ for which the best value of the PSNR is reached 
when the image is processed using the SK algorithm;
\item   "PSNR - SK"  denotes the values of the PSNR achieved by the image processed by the SK algorithm for 
$(N,w)_{Psnr}^{\max}$.
\end{itemize}
Note that, in all the above cases the PSNR  is computed using 
 the original  image of dimension $N \times M$  as the reference image.
{\footnotesize \phantom{a} \hskip-2cm 
\begin{table}[!h]
\begin{center}
\caption{ The numerical values of the PSNR. } \label{tab1}
\begin{tabular}{l c c | c | c c  }
\toprule  
{\bf Image}  &  
$\delta_{PSNR}^{max}$ &
 \specialcell{\bf PSNR \\ $\delta_{Psnr}^{\max}$ \\ \bf Fuzzy}
 &   \specialcell{\bf PSNR \\ \bf bicubic \\ \phantom{a}} &
\specialcell{$(N,w)_{Psnr}^{\max}$ \\  } 
& \specialcell{\bf PSNR\\ $(N,w)_{Psnr}^{\max}$\\ \bf SK }

    \\
\midrule
{\bf Baboon}  &  0.86 & 22.0121 &   22.2143  & (2,10) (see Table \ref{bsk-psnr}) & 22.43814 \\ 
\midrule
{\bf Boat}       &  0.72 & 29.6480 &  28.0849  & (4,20) (see Table \ref{boat-psnr}) & 29.3667\\ 
\midrule
{\bf City}        &  0.58 & 24.9440 &  24.4878 & (4,20) (see Table \ref{city-psnr}) & 25.08078  \\ 
\midrule
{\bf Mountain}  &  0.55 & 22.4655 & 21.4886  & (6,20) (see Table \ref{mountain-psnr}) & 23.35247   \\ 
\bottomrule
\end{tabular}
\end{center}
\end{table}
}
%

 Moreover, in Table  \ref{tab2}, which has
 the same meaning as in Table \ref{tab1},   the values of the likelihood index S are listed and analysed,
 for  the case of the fuzzy-type,  bicubic  and   SK algorithms.\\
%
{\footnotesize
\begin{table}[!h]
\caption{ The numerical values of the likelihood index S.} \label{tab2}
\begin{center}
\begin{tabular}{l  c  c |  c | c  c   }
\toprule  
{\bf Image}   & $\delta_{S}^{\max}$     & \specialcell{\bf S index \\ $\delta_S^{\max}$  \\  \bf Fuzzy } 
		         &   \specialcell{ \bf S index \\ bicubic\\  \phantom{a}} 
 	              & \specialcell{$(N,w)_{S}^{\max}$ \\ }
	       	  & \specialcell{\bf S index\\  $(N,w)_{S}^{\max}$ \\ \bf SK } 
		    \\ \midrule
{\bf Baboon} &  0.86 & 0.9456 &   0.9449  & (7,25) (see Table \ref{bsk-sindex}) & 0.947694
 \\  \midrule 
{\bf Boat}  &  0.78 & 0.98051 &   0.9771  & (4,20) (see Table \ref{boat-sindex}) & 0.980565
\\  \midrule
{\bf City } &  0.58 & 0.9735 &   0.9711  & (4,20) (see Table \ref{city-sindex}) & 0.974093
 \\  \midrule
{\bf Mountain } &  0.61 & 0.9523 &   0.9455  & (5,20) (see Table \ref{mountain-sindex}) & 0.955541
 \\  \bottomrule
\end{tabular}
\end{center}
\end{table}
}

%
 In Tables \ref{tab1} and \ref{tab2}   we observe a  similar trend in performances with respect to the two indices 
with the exception of the boat image in which the fuzzy algorithm performs better than the SK at least with respect to the PSNR 
index and in  the  baboon image where the best resolutions for the SK algorithm are  obtained for  $(N,w)_{Psnr}^{\max}$ 
and $(N,w)_{S}^{\max}$ which are very distant from each other.\\

Finally,  an analysis concerning the CPU time employed by each  of the considered algorithms to process any single image can be performed. 
 The CPU times are listed in Table \ref{tab3}.
Since it is not the purpose of the present study to determine all the  CPU times,  we consider and  compare only  the times of the best reconstructions.
Therefore in Table \ref{tab3} the values of CPU times are considered for the reconstructed images obtained for 
$\delta_{Psnr}^{\max}$, $\delta_{S}^{\max}$, for the fuzzy algorithm and  
$(N,w)_{Psnr}^{\max}$, $(N,w)_{S}^{\max}$, for the SK algorithm
and quoted in Tables \ref{tab2} and \ref{tab3}.

\begin{table}[!h]
\footnotesize
    \caption{ The CPU for the rescaled images of the best approximations for   the PSNR and  the likelihood index S.}
    \label{tab3}
\phantom{a} \hskip-1.62cm
    \begin{tabular}{ l *{6}{>{\raggedleft\arraybackslash}p{1.62cm}} }
    \toprule
        \\ 
        &
        \multicolumn{6}{c}{\textbf{\specialcell{The CPU time for the best approximations \\ with respect to PSNR and S indices} }} \\ 
    \cmidrule{2-7} 
    \multirow{-2}{*}{\parbox{1.45cm}{\textbf{Case n=3}}} &
        \multicolumn{1}{c}{\textbf{\specialcell{dim \\ image}}} &
        \multicolumn{1}{c}{\textbf{bicubic}} &
      \multicolumn{1}{c}{\textbf{\specialcell{Fuzzy\\$\delta_{Psnr}^{\max}$}  }} &
        \multicolumn{1}{c}{\textbf{\specialcell{SK\\ $(N,w)_{Psnr}^{\max}$}}} &
        \multicolumn{1}{c}{\textbf{\specialcell{Fuzzy\\ $\delta_{S}^{\max}$}}} &
        \multicolumn{1}{c}{\textbf{\specialcell{SK\\ $(N,w)_{S}^{\max}$}}} 
 \\  \midrule
    \textbf{Baboon}   &  255 $\times$ 255 &  0.054889  &  {2.532010}  &  192.986372  &    2.532010 &    5.350619  \\ \midrule
    \textbf{Boat}       & 504 $\times$ 504  &   0.091244  &  0.694429 &  40.098191  &  0.585434 &    40.098191  \\ \midrule
 \textbf{City}           & 675 $\times$ 900  &   0.100341   &  1.137932 &  192.654707 &  1.137932 &    192.654707  \\ \midrule
    \textbf{Mountain} &  450 $\times$ 600  &   0,084561 &  0.690313 &  22.639233 &  0.595628 &    17.945265  \\ 
      \bottomrule
    \end{tabular}
\end{table}


\begin{remark} \rm
Note that as
 the parameter  $N$  increases, the order of decay of the Jackson kernel increase as well 
and therefore the CPU  time of the SK algorithm
 decreases;  however, the similarity indices worsen.
The SK algorithm is the most expensive
 from the point of view of CPU time. If, however, instead of considering its best approximation, we take into account the values of $(N,k)$
 so that $N$ is large enough and the SK algorithm performs better than the fuzzy one, we can strongly reduce the 
 CPU time. 
For example
 if we consider the reconstructed image of Baboon with $N=12$ and $w=15$ we need $2.857455s$ and we obtain, 
 accordingly to Table \ref{bsk-psnr}, a better result with respect to the fuzzy algorithm in a much shorter time than   $(N,w)_{Psnr}^{\max}$.
However if we look at Tables \ref{bsk-psnr}-\ref{mountain-sindex}, in the Appendix,we can see that, except for the boat,
 the values of the PSNR or S-index  for the SK algorithm are better than those of the fuzzy algorithm when
 $N$ is quite large, so, in this case, the CPU times of the SK algorithm decrease, continuing to achieve better performances.
\end{remark}

For the sake of completeness, we also considered the case of the application of the  abovementioned rescaling algorithms by a
 resize factor   $R = (2k+1), \,  k=2, 3$.\\
In practice, we repeated the above experiments reducing the considered original images to 
 have  dimension  of  
 ${N \over {2k+1}} \times {M \over {2k+1}},\,  k=2,3$. Consequently, 
by the above  methods,   they have been processed in order to reobtain  images scaled to the original dimension.\\ 

Due to the compatibility between the amplitude of the scale factor and the dimensions of the original images, in this case we 
   considered   the images
  "baboon" and "mountain" for $k=2$ and "boat" for $k=3$. 
For the application of the SK algorithm.\\ 

The corresponding numerical results of this case  are
 presented in Figures \ref{fig2-r5} and \ref{fig2-r7}. 
 Additionally here, the fuzzy-type algorithm is applied for every $\delta$ between $0$ and $1$, with the same step-size of $0.01$.\\
%

{
\begin{figure}[h!]
\fbox{\includegraphics[scale=0.262]{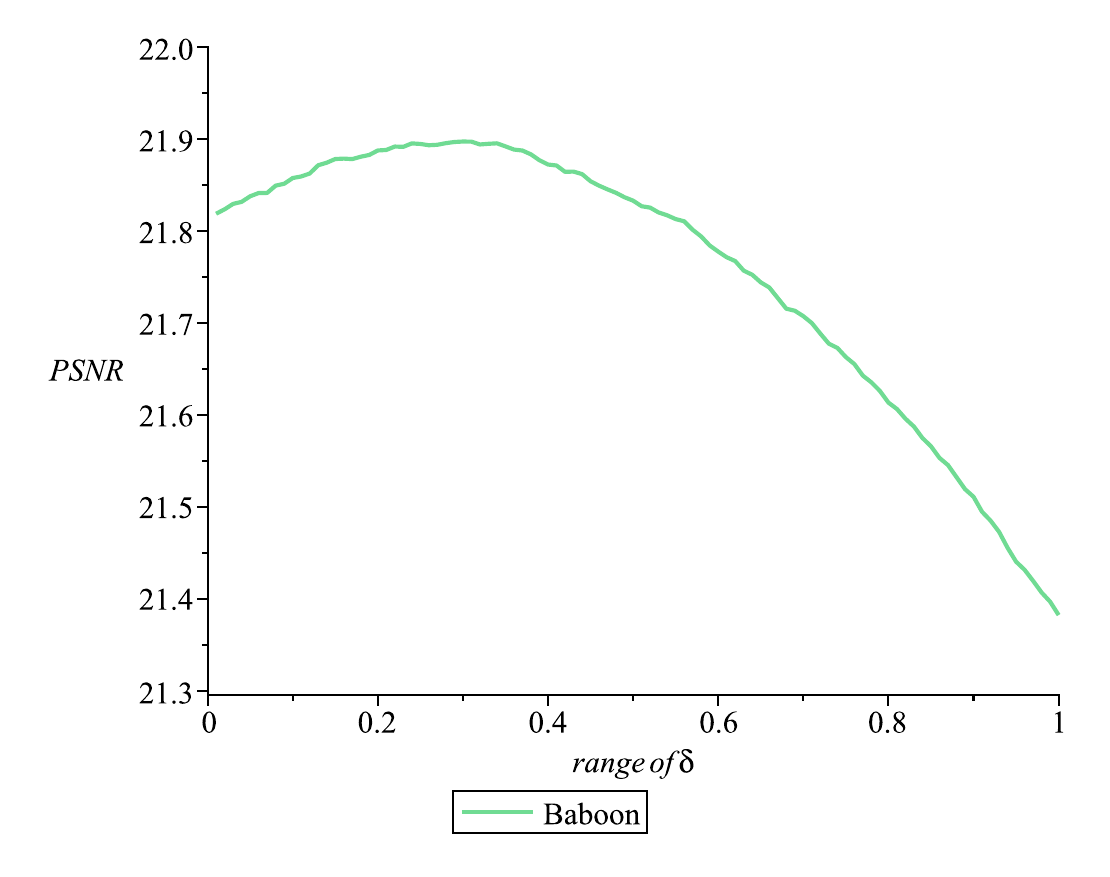}
\hskip0.3cm
\includegraphics[scale=0.262]{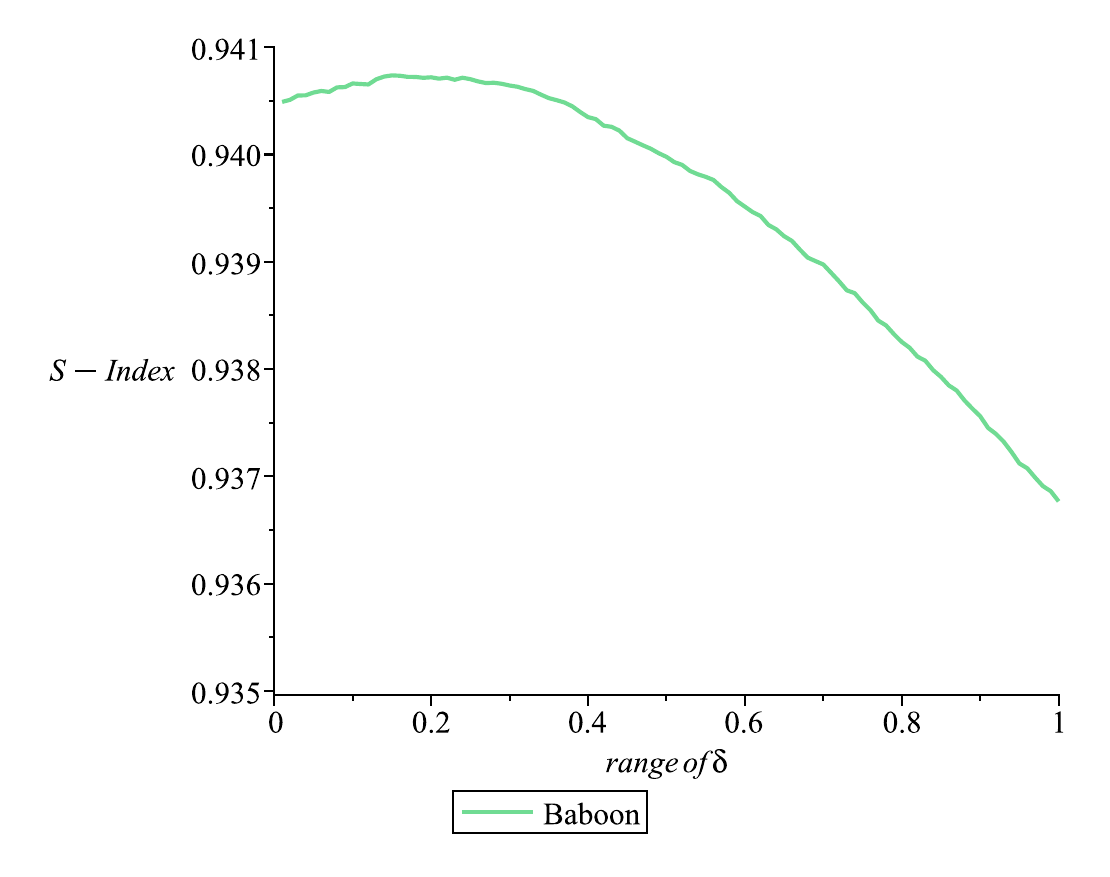} \quad }
\vskip0.4cm
\fbox{\includegraphics[scale=0.252]{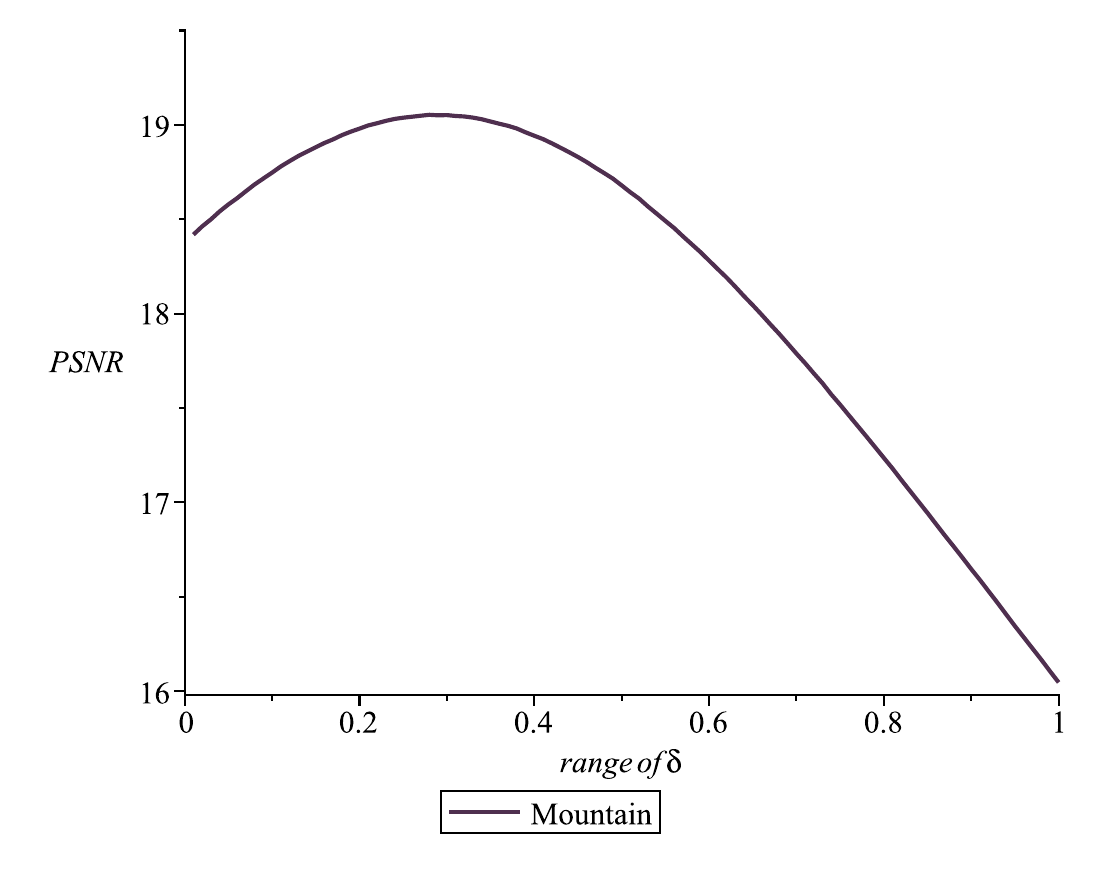}
\hskip0.3cm
\includegraphics[scale=0.27]{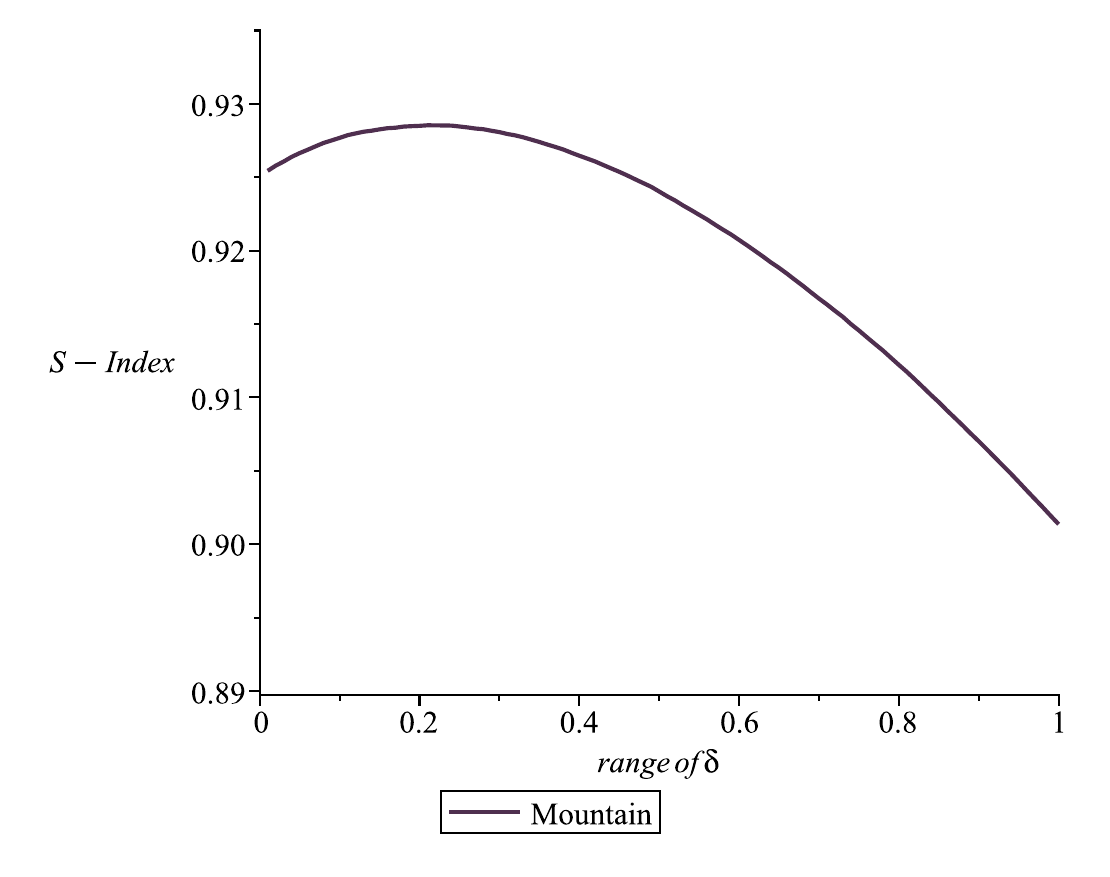} \quad}
\caption{\footnotesize $R=5$: the plots of the values of the PSNR and likelihood index S computed for  of the  reconstructed images of
  baboon and mountain with  the fuzzy-type 
algorithm when the parameter $\delta$  varies from $0$ to $1$ with step-size of $0.01$.} \label{fig2-r5}
\end{figure}
}

{
\begin{figure}[h!]\label{delta-fig7}
\fbox{\includegraphics[scale=0.27]{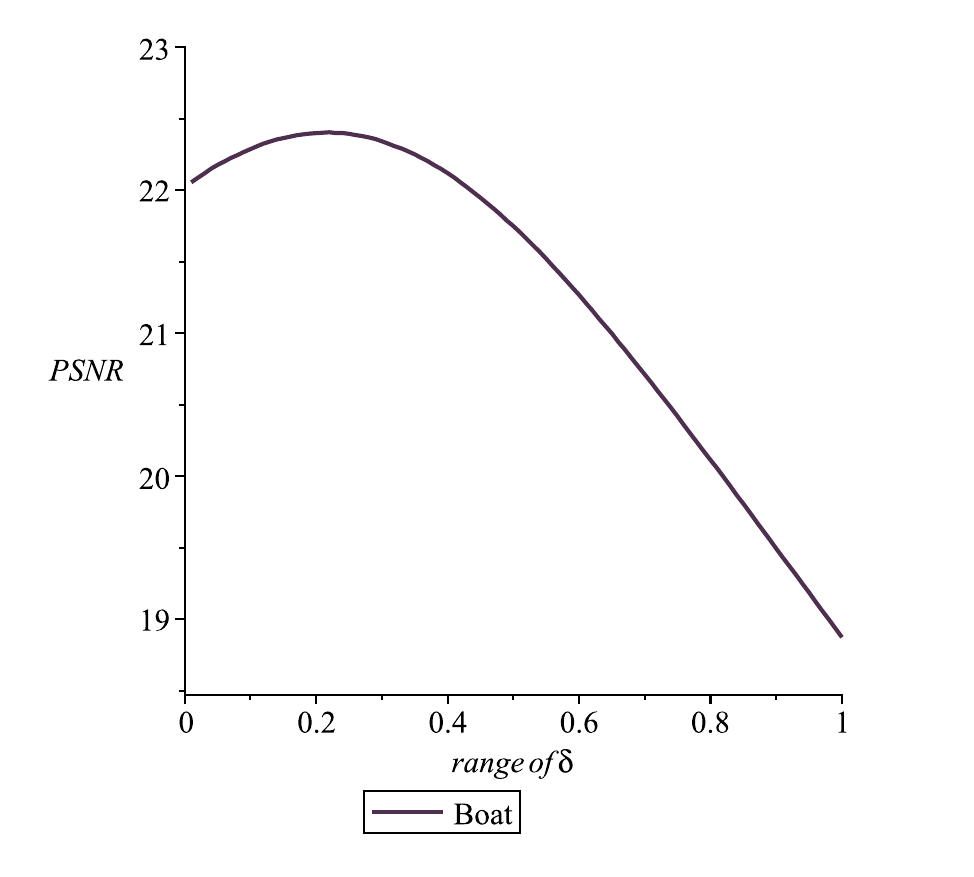}
\hskip0.35cm
\includegraphics[scale=0.27]{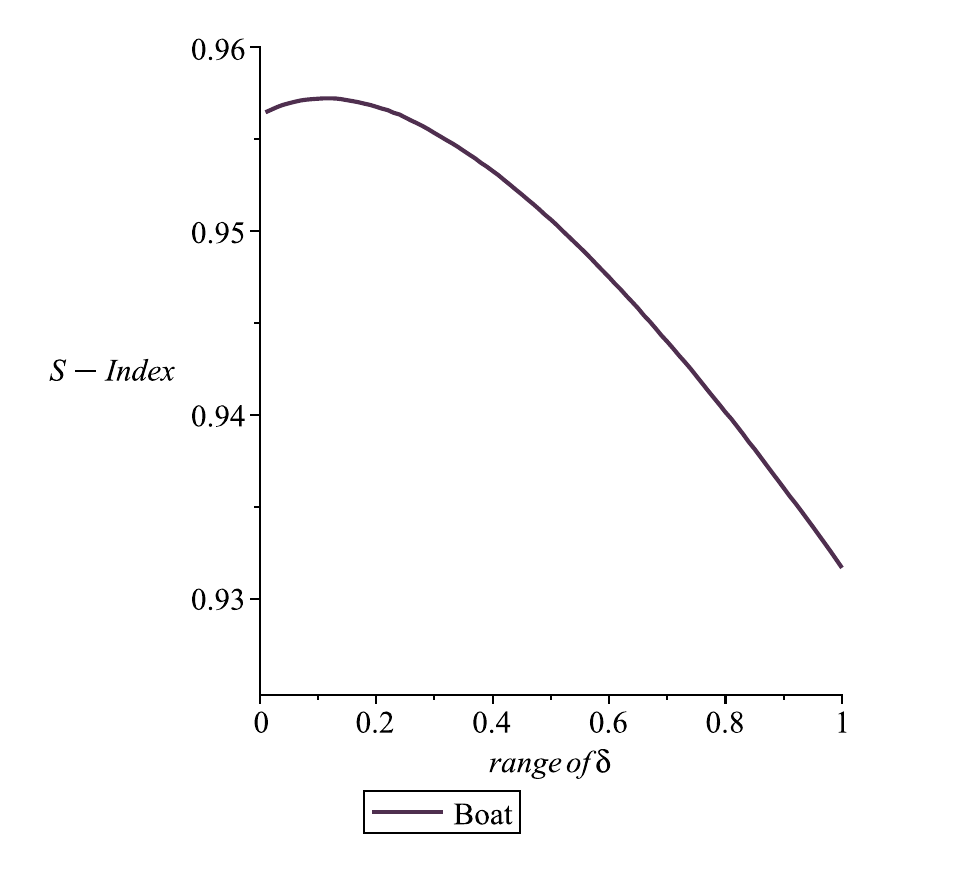} \quad }
\caption{\footnotesize $R=7$: the plots of the values of the PSNR and likelihood index S computed for the  reconstructed images
  of boat with  the fuzzy-type 
algorithm when the parameter $\delta$  varies from $0$ to $1$ with step-size of $0.01$.} \label{fig2-r7}
\end{figure}
}
\newpage

\begin{remark} \rm 
The performances of the fuzzy-type algorithm are dependent  on the value of the parameter $\delta$. In all the considered cases,
 it seems that the curves of the PSNR and S index plots are both concave and achieve a maximum approximatively in the middle 
zone of the interval $[0,1]$, if we consider the experiments  with a scaling factor equal to 3.
 This fact  seems to be more evident in the figures: 
 for Boat, City, and Mountain. 
 When the scaling factor is equal to 5 or 7, the point of the maximum shifts toward the left, as shown in the following Figures \ref{fig2-r5-delta} and \ref{fig2-r7-delta}. 
{
\begin{figure}[h!]
\quad \fbox{\includegraphics[scale=0.25]{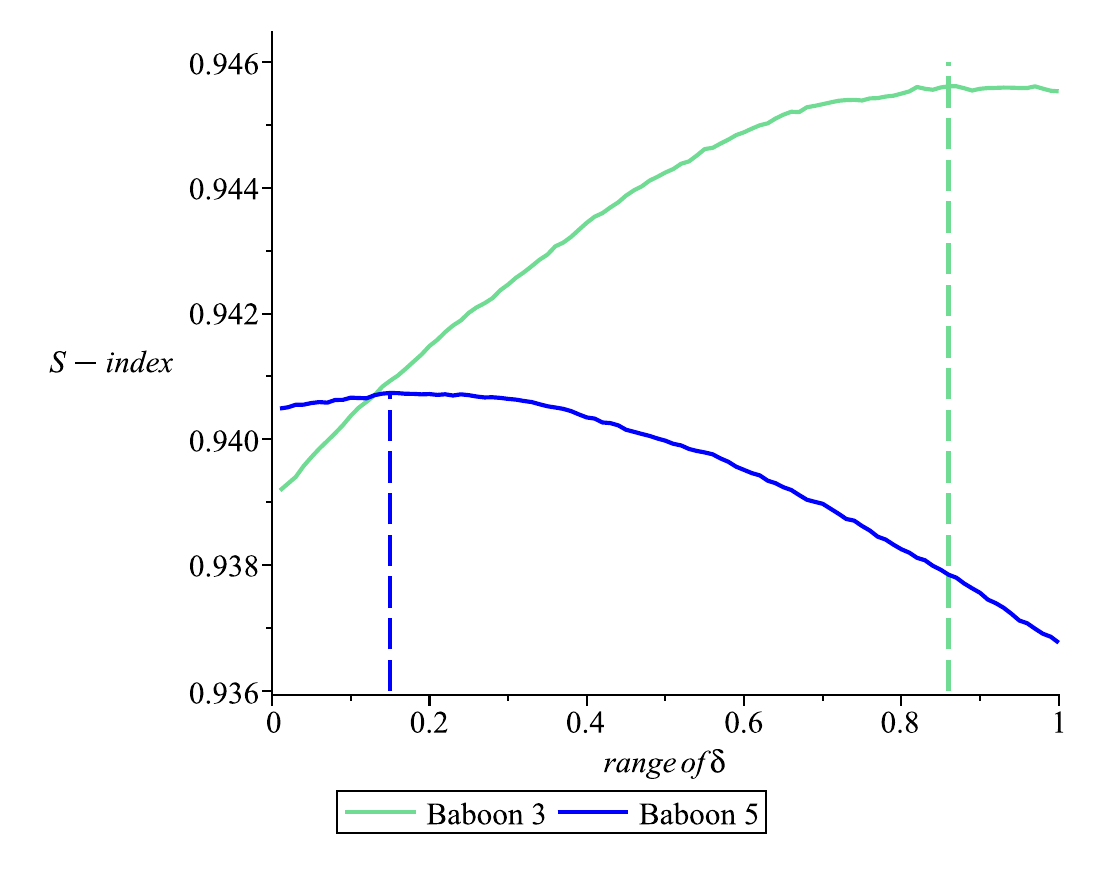}\,
\hskip0.3cm
\includegraphics[scale=0.25]{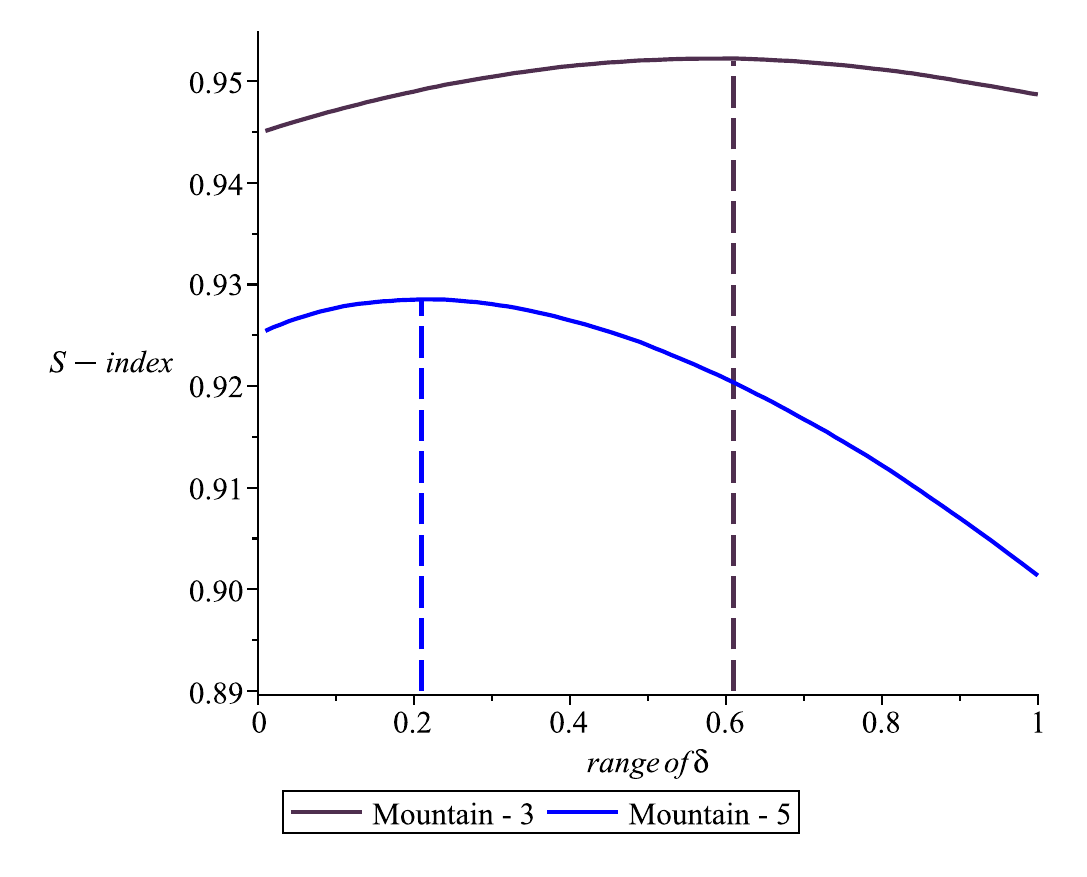}  }
\caption{\footnotesize $R=3,5$: the left shift of the  maximum in  the plots of the values of the likelihood index S computed for the  reconstructed images  of the baboon and the mountain  with  the fuzzy-type 
algorithm when the parameter $\delta$  varies from $0$ to $1$ with step-size of $0.01$.} \label{fig2-r5-delta}
\end{figure}
}
\vskip-.3cm
{
\begin{figure}[h!]
\quad \fbox{\includegraphics[scale=0.22]{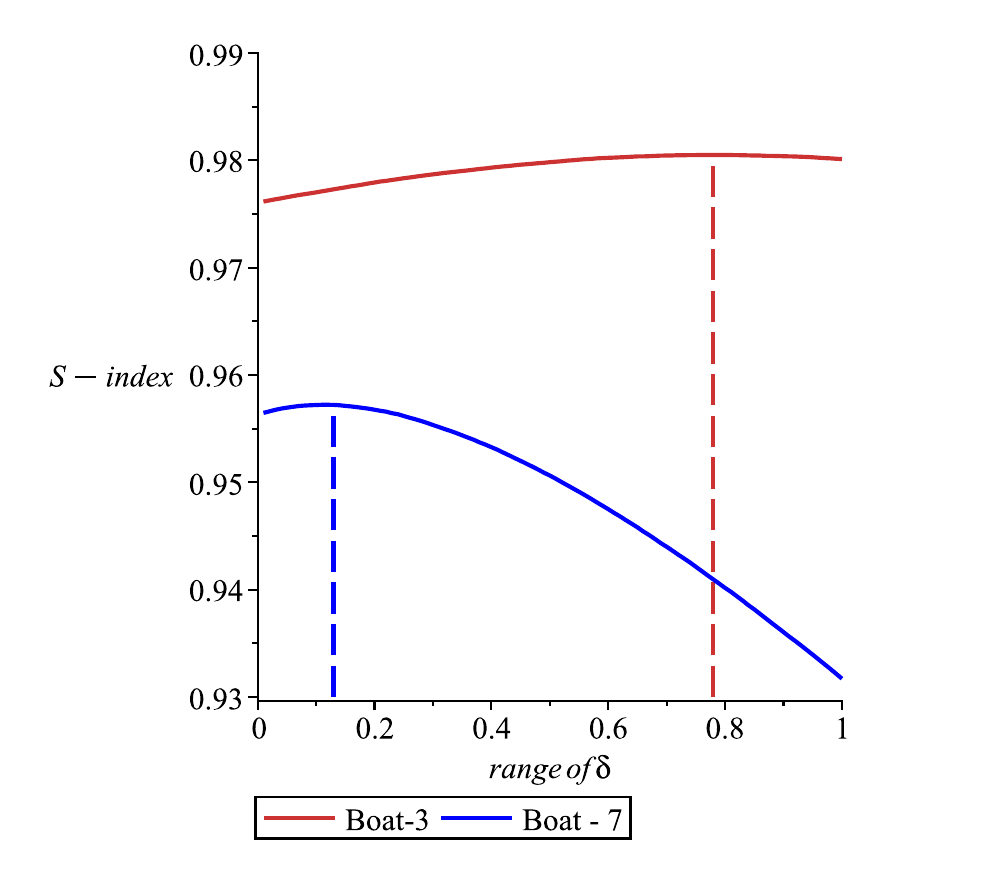} \, }
\caption{\footnotesize  $R=3,7$: the left shift of the  maximum in  the plots of the values of the likelihood index S computed for the 
 reconstructed images  of the boat  with  the fuzzy-type 
algorithm when the parameter $\delta$  varies from $0$ to $1$ with step-size of $0.01$.} \label{fig2-r7-delta}
\end{figure}
}
\end{remark}


{
\begin{figure}[h!]
\fbox{\includegraphics[scale=0.4]{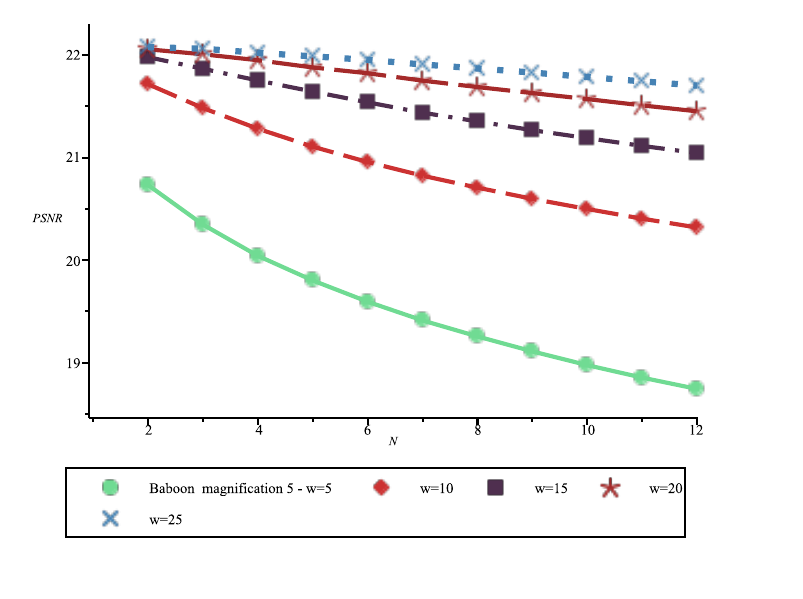} \, 
 \includegraphics[scale=0.4]{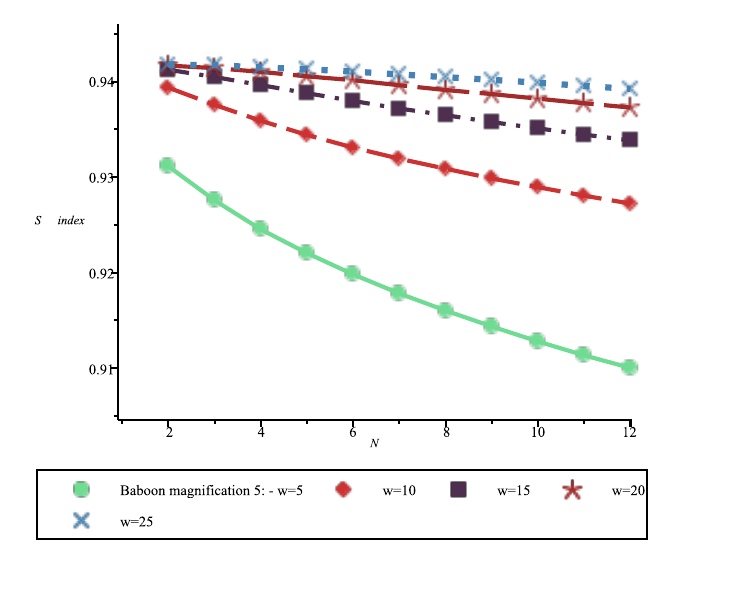} \,\,}
\\ \vskip.42cm
\fbox{\includegraphics[scale=0.3]{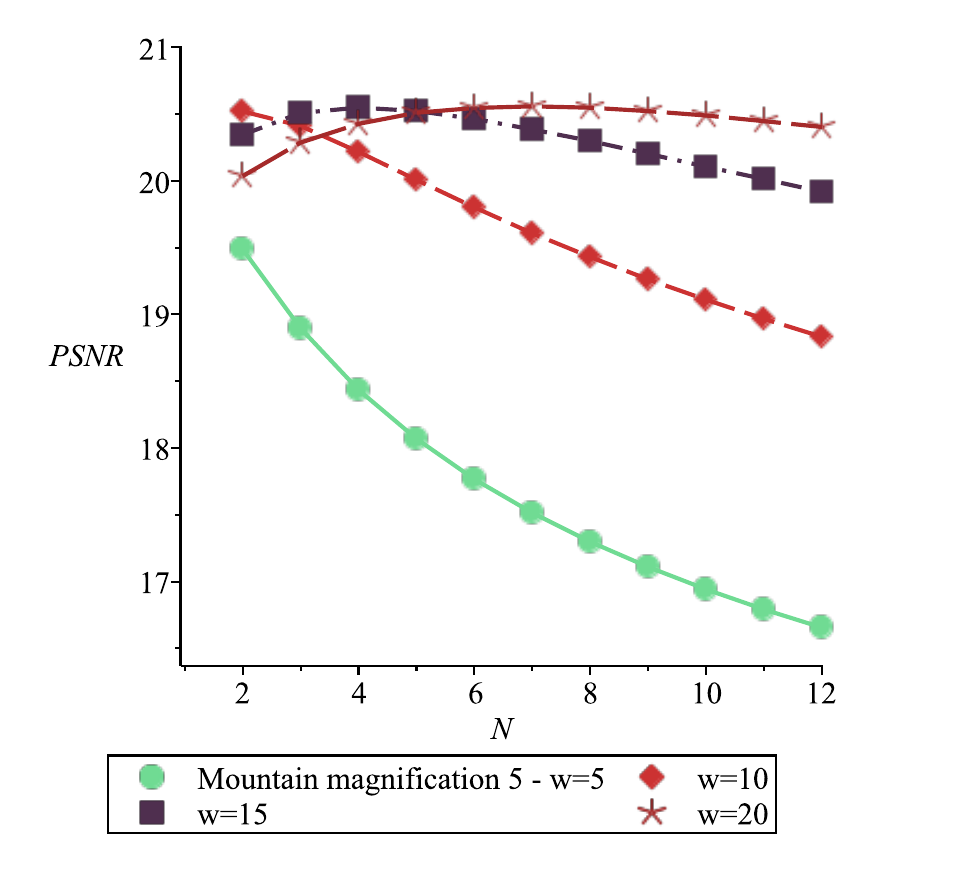} \,\,   
\includegraphics[scale=0.3]{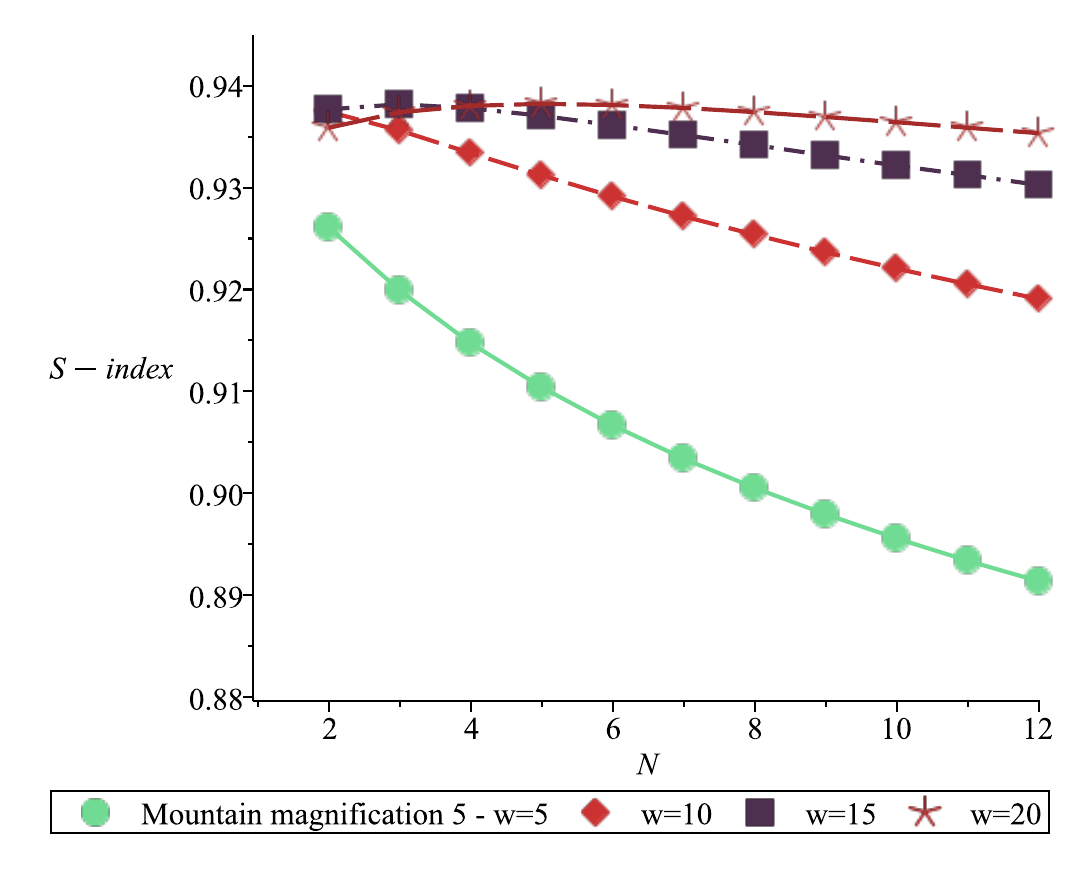} \,\,}
\caption{\footnotesize The plots of the values of the computed indices for the  reconstructed baboon and mountain  images when the magnification is 5} \label{fig-baboon5}
\end{figure}
}

\begin{remark}\rm
From the plots of Figure \ref{fig-baboon5} and the data of Tables \ref{bsk-psnr-5} - \ref{msk-sindex-5} it seeems clear that globally the values of the considered indices are smaller than those achieved in the
 corresponding cases with magnification factor equal to 3.\\
This seems quite natural since, we are starting the reconstructions from images that are sensibily smaller than those used in the previous reconstructions, 
and this can be translated into a process that is based on much less starting information, with respect to the previous case,  that can difficulty produce accurate results.
Following this reasoning, we can also justify the fact that increasing the value of $w$ the quality of the reconstruction does not improve.
\end{remark}

\begin{table}[h!]
\tiny
    \caption{\footnotesize The numerical values of the PSNR and  the likelihood index S, when the images are rescaled by a factor equal to 5.
 The  values must be interpreted as in the previous tables.}
    \label{tab4a}
\hspace*{-1.7cm} \tiny
    \begin{tabular}{ l *{8}{>{\raggedleft\arraybackslash}p{1.2cm}} }    \toprule
        &
        \multicolumn{8}{c}{\textbf{The numerical values of the PSNR and  the likelihood index S}} \\ 
 \toprule
        &     \multicolumn{4}{c}{\textbf{PSNR index}} &    \multicolumn{4}{c}{\textbf{S index}} \\       
    \cmidrule(r){2-5} \cmidrule(l){6-9}    \\    
    \multirow{-2}{*}{\parbox{1cm}{\textbf{Case n=5}}} &
        \multicolumn{1}{c}{\textbf{$\delta^{\max}$}} &
        \multicolumn{1}{c}{\textbf{\specialcell{Fuzzy\\ $\delta^{\max}$ }}} &
        \multicolumn{1}{c}{\textbf{\specialcell{ $(N,w)_{\max}$}}} &
       \multicolumn{1}{c}{\textbf{\specialcell{SK}}} &
        \multicolumn{1}{c}{\textbf{$\delta^{\max}$}} &
        \multicolumn{1}{c}{\textbf{\specialcell{Fuzzy\\ $\delta^{\max}$ }}} &
        \multicolumn{1}{c}{\textbf{\specialcell{$(N,w)_{\max}$}}} & 
           \multicolumn{1}{c}{\textbf{\specialcell{SK\\ }}} \\
 \cmidrule(r){2-5} \cmidrule(r){6-9}
    \textbf{Baboon}     &  0.30   & 21.897  & (2,25)  &  22.0775 &   0.15    & 0.9407  & (2,25)  &  0.94178 \\
    \cmidrule(r){2-5} \cmidrule(l){6-9}    \\    
    \textbf{Mountain}   &0.30  & 19.05169  & $(7,20)$ & 20.55734   &  0.21  &  0.92854  & $(5,20)$ & 0.93822   \\    
      \bottomrule
    \end{tabular}
\end{table}


\section{Appendix}\label{Appendix}

\subsection{magnification 3}
\phantom{a} \\
We report here the values of the two similarity indices for  the rescaled images (magnification  3) using the SK algorithm depending 
on the values of the two parameters $N$ and $w$. Each table  refers to a single image and a single index of similarity.
The values of $(N,w)_{PSNR}$ and $(N,w)_S$ for which the maximum of the similarity indices is reached  appear in the tables in bold.

 \begin{table}[htp]
\small
\caption{ \footnotesize The numerical values of the PSNR for the baboon image}\label{bsk-psnr}
\begin{tabular}{rccccc}
\toprule
N  & w=5      & {\bf w=10 }    & w=15     & w=20     & w=25     \\ \midrule
{\bf 2}  & 21,85016 & {\bf 22,43814 } & 22,32794 & 22,10527 & 21,83349 \\ \midrule
3  & 21,48853 & 22,37021 & 22,41231 & 22,28217 & 22,08483 \\ \midrule
4  & 21,20701 & 22,26074 & 22,43747 & 22,36078 & 22,2281  \\ \midrule 
5  & 20,97893 & 22,14198 & 22,42497 & 22,40758 & 22,30682 \\ \midrule
6  & 20,78909 & 22,02888 & 22,39183 & 22,4312  & 22,35847 \\ \midrule 
7  & 20,61677 & 21,91724 & 22,35921 & 22,43585 & 22,38839 \\ \midrule
8  & 20,47379 & 21,81542 & 22,29998 & 22,43023 & 22,41282 \\ \midrule
9  & 20,34389 & 21,7138  & 22,24579 & 22,41857 & 22,42535 \\ \midrule 
10 & 20,22084 & 21,62172 & 22,19483 & 22,40215 & 22,43644 \\ \midrule 
11 & 20,11594 & 21,5353  & 22,13739 & 22,38263 & 22,43331 \\ \midrule
12 & 20,01313 & 21,452   & 22,09212 & 22,36252 & 22,43133 \\ \bottomrule
\end{tabular}
\end{table}

\begin{table}[h!]
\small
\caption{ \footnotesize  The numerical values of the S-index for the baboon image}\label{bsk-sindex}
\begin{tabular}{rccccc}
\toprule
\multicolumn{1}{c}{N} & w=5      & w=10     & w=15     & {\bf w=20}     & w=25     \\ \midrule
2                     & 0,940944 & 0,947122 & 0,947393 & 0,94603  & 0,944384 \\ \midrule
3                     & 0,937897 & 0,945927 & 0,947648 & 0,947146 & 0,945926 \\ \midrule
4                     & 0,935421 & 0,944679 & 0,947305 & 0,94758  & 0,946816 \\ \midrule
{\bf 5}                    & 0,933376 & 0,943534 & 0,946765 &  0,947678 & 0,947299 \\ \midrule
6                     & 0,931641 & 0,942475 & 0,946223 & 0,947527 & 0,947564 \\ \midrule
7                     & 0,930041 & 0,94149  & 0,945736 & 0,94727  &{\bf 0,947694} \\ \midrule
8                     & 0,928687 & 0,940611 & 0,945117 & 0,946969 & 0,947649 \\ \midrule
9                     & 0,927407 & 0,939758 & 0,944554 & 0,946683 & 0,947556 \\ \midrule
10                    & 0,926217 & 0,938982 & 0,944042 & 0,94636  & 0,947444 \\ \midrule
11                    & 0,925169 & 0,938254 & 0,9435   & 0,946073 & 0,947228 \\ \midrule
12                    & 0,924114 & 0,937539 & 0,943055 & 0,945802 & 0,947034 \\\bottomrule
\\
\end{tabular}
\end{table}


\begin{table}[h!]
\small
\caption{ \footnotesize  The numerical values of the PSNR for the boat image}\label{boat-psnr}
\begin{tabular}{rcccc}
\toprule
N  & w=5      & w=10     & w=15     & {\bf w=20}       \\ \midrule
2  & 26,1264  & 28,88491 & 29,33233 & 29,1455    \\ \midrule
3  & 25,10887 & 28,2995  & 29,25153 & 29,35566   \\ \midrule
{\bf 4 } & 24,39938 & 27,69831 & 29,01746 & {\bf 29,3667}    \\ \midrule
5  & 23,8721  & 27,18115 & 28,73699 & 29,29081   \\ \midrule
6  & 23,45101 & 26,74578 & 28,4665  & 29,16879   \\ \midrule
7  & 23,11319 & 26,36313 & 28,18154 & 29,02029  \\ \midrule
8  & 22,82393 & 26,03003 & 27,90549 & 28,86864   \\ \midrule
9  & 22,57538 & 25,73335 & 27,65742 & 28,70946   \\ \midrule
10 & 22,35332 & 25,46989 & 27,4174  & 28,5503    \\ \midrule
11 & 22,15799 & 25,23723 & 27,20233 & 28,39299   \\ \midrule
12 & 21,98237 & 25,02072 & 26,99476 & 28,23009 
       \\ \bottomrule
& & & &
\end{tabular}
\end{table}
\phantom{a} \\

\begin{table}[h!]
\small
\caption{ \footnotesize  The numerical values of the S-index for the boat image}\label{boat-sindex}
\begin{tabular}{lllll}
\toprule
N  & w=5      & w=10     & w=15     & {\bf w=20}     \\ \midrule
2  & 0,971368 & 0,979557 & 0,980425 & 0,979874 \\ \midrule
3  & 0,967317 & 0,978083 & 0,980413 & 0,98041  \\ \midrule
{\bf 4}  & 0,964053 & 0,976406 & 0,979877 & {\bf 0,980565} \\ \midrule
5  & 0,961385 & 0,974849 & 0,979184 & 0,980472 \\ \midrule
6  & 0,95908  & 0,973451 & 0,978483 & 0,980232 \\ \midrule
7  & 0,957091 & 0,972149 & 0,97776  & 0,979876 \\ \midrule
8  & 0,955326 & 0,970935 & 0,97699  & 0,97951  \\ \midrule
9  & 0,953716 & 0,969804 & 0,976281 & 0,979102 \\ \midrule
10 & 0,952239 & 0,968759 & 0,975565 & 0,978695 \\ \midrule
11 & 0,950878 & 0,96779  & 0,97491  & 0,978308 \\ \midrule
12 & 0,949631 & 0,966857 & 0,974246 & 0,97789 \\ \bottomrule
\\
\end{tabular}
\end{table}
\phantom{a}
\clearpage
\begin{table}[h!]
\small
\caption{ \footnotesize  The numerical values of the PSNR for the city image}\label{city-psnr}
\begin{tabular}{rcccc}
\toprule
N  & w=5      & w=10     & w=15     & {\bf w=20 }    \\ \midrule
2  & 23,54893 & 24,94082 & 25,04948 & 24,86407 \\ \midrule
3  & 22,9495  & 24,68724 & 25,06223 & 25,03318 \\ \midrule
{\bf 4}  & 22,51778 & 24,40862 & 24,99417 & {\bf 25,08078} \\ \midrule
5  & 22,18384 & 24,1397  & 24,89032 & 25,07869 \\ \midrule
6  & 21,91455 & 23,89762 & 24,76533 & 25,0469  \\ \midrule
7  & 21,68927 & 23,68611 & 24,64658 & 25,00216 \\ \midrule
8  & 21,49285 & 23,50061 & 24,51943 & 24,94458 \\ \midrule
9  & 21,31981 & 23,32962 & 24,39139 & 24,87755 \\ \midrule
10 & 21,16648 & 23,17378 & 24,26777 & 24,80927 \\ \midrule
11 & 21,02844 & 23,03034 & 24,15284 & 24,74202 \\ \midrule
12 & 20,90149 & 22,90277 & 24,03818 & 24,67353 \\ \bottomrule
& & & &
\end{tabular}
\end{table}
\phantom{a} \\  \vskip.2cm
\begin{table}[h!]
\small
\caption{ \footnotesize  The numerical values of the S-index for the city image}\label{city-sindex}
\begin{tabular}{rcccc}
\toprule
N  & w=5      & w=10     & w=15     & {\bf w=20}     \\ \midrule
2  & 0,966276 & 0,973093 & 0,974044 & 0,973677 \\ \midrule
3  & 0,963224 & 0,971649 & 0,973883 & 0,974052 \\ \midrule
{\bf 4}  & 0,96082  & 0,970291 & 0,973345 &{\bf  0,974093} \\ \midrule
5  & 0,958814 & 0,969054 & 0,972696 & 0,973956 \\ \midrule
6  & 0,957099 & 0,967914 & 0,972007 & 0,973693 \\ \midrule
7  & 0,955605 & 0,966908 & 0,971388 & 0,973352 \\ \midrule
8  & 0,95425  & 0,965991 & 0,970786 & 0,972984 \\ \midrule
9  & 0,953023 & 0,965129 & 0,970177 & 0,972603 \\ \midrule
10 & 0,951876 & 0,964339 & 0,9696   & 0,972217 \\ \midrule
11 & 0,950819 & 0,963593 & 0,969084 & 0,971851 \\\midrule
12 & 0,949828 & 0,962915 & 0,968562 & 0,971497 \\ \bottomrule
\\
\end{tabular}
\end{table}
\clearpage

\begin{table}[h!]
\small
\caption{ \footnotesize  The numerical values of the PSNR for the mountain image}\label{mountain-psnr}
\begin{tabular}{rcccc}
\toprule 
N  & w=5      & w=10        & w=15     & {\bf w=20}     \\ \midrule
2  & 21,76631 & 23,27749 & 23,1204  & 22,68569 \\ \midrule
3  & 21,03285 & 23,07359 & 23,31022 & 23,05427 \\ \midrule
4  & 20,47727 & 22,76765 & 23,33451 & 23,23154 \\ \midrule
5  & 20,04889 & 22,47274 & 23,2664  & 23,31831 \\ \midrule
{\bf 6}  & 19,69711 & 22,19207 & 23,1556  & {\bf 23,35247} \\ \midrule
7  & 19,40644 & 21,93179 & 23,03366 & 23,34492 \\ \midrule
8  & 19,15224 & 21,70653 & 22,89272 & 23,31629 \\ \midrule
9  & 18,93389 & 21,50075 & 22,75439 & 23,26504 \\ \midrule
10 & 18,7375  & 21,30299 & 22,61814 & 23,19806 \\ \midrule
11 & 18,55898 & 21,12716 & 22,4842  & 23,13559 \\ \midrule
12 & 18,40131 & 20,96599 & 22,35643 & 23,06147 \\ \bottomrule
\end{tabular}
\end{table}


\begin{table}[h!]
\small
\caption{ \footnotesize The numerical values of the S-index for the mountain image}\label{mountain-sindex}
\begin{tabular}{lllll}
\toprule
N  & w=5      & w=10        & w=15     & {\bf w=20 }    \\  \midrule
2  & 0,943975 & 0,954693 & 0,955022 & 0,953237 \\  \midrule
3  & 0,938593 & 0,952793 & 0,955499 & 0,954773 \\  \midrule
4  & 0,934135 & 0,950659 & 0,95508  & 0,955388 \\  \midrule
{\bf 5}  & 0,930443 & 0,948655 & 0,954247 & {\bf 0,955541} \\  \midrule
6  & 0,927252 & 0,946779 & 0,953304 & 0,955412 \\  \midrule
7  & 0,92449  & 0,945008 & 0,952377 & 0,955075 \\  \midrule
8  & 0,921978 & 0,943414 & 0,951415 & 0,954651 \\  \midrule
9  & 0,919735 & 0,941943 & 0,950494 & 0,954142 \\  \midrule
10 & 0,917661 & 0,940508 & 0,949581 & 0,953597 \\  \midrule
11 & 0,915726 & 0,939196 & 0,948674 & 0,953094 \\  \midrule
12 & 0,913969 & 0,93795  & 0,947829 & 0,952548\\ \bottomrule
\end{tabular}
\end{table}

\clearpage
{
\begin{table}[h!]
\footnotesize
\caption{ \footnotesize The numerical values of the S-index for the boat image.
}
\label{boat-num}
\hspace*{-.3cm}
\begin{tabular}{c c c c}
\toprule
 {[}0.01; 0.97616558969266853{]} & {[}0.02; 0.97626759036948307{]} & {[}0.03; 0.97637292571349765{]} & {[}0.04; 0.97646385714780093{]} \\ \midrule
 {[}0.05; 0.97656766410268026{]} &
{[}0.06; 0.97666971109431622{]} & {[}0.07; 0.97676601504799732{]} & {[}0.08; 0.97684976768485399{]} \\ \midrule
 {[}0.09; 0.97693415329094790{]} & {[}0.10; 0.97701875503287983{]} &
{[}0.11; 0.97711845540685693{]} & {[}0.12; 0.97720917070532198{]} \\ \midrule
 {[}0.13; 0.97730444029463814{]} & {[}0.14; 0.97739659135259349{]} & {[}0.15; 0.97748522248405789{]} &
{[}0.16; 0.97758537056798234{]} \\ \midrule
 {[}0.17; 0.97766006093795910{]} & {[}0.18; 0.97775185691560795{]} & {[}0.19; 0.97785212850572534{]} & {[}0.20; 0.97793807337750072{]} \\ \midrule
{[}0.21; 0.97803451627564830{]} & {[}0.22; 0.97809214735277061{]} & {[}0.23 ; 0.97817858624931775{]} & {[}0.24; 0.97825302960690674{]} \\ \midrule
 {[}0.25; 0.97833929868243896{]} &
{[}0.26; 0.97840464889659817{]} & {[}0.27; 0.97848502055143227{]} & {[}0.28; 0.97855926321145992{]} \\ \midrule
 {[}0.29; 0.97862980068570926{]} & {[}0.30; 0.97869913397458064{]} &
{[}0.31; 0.97876408279361438{]} & {[}0.32; 0.97883939069455217{]} \\ \midrule
 {[}0.33; 0.97889325483279954{]} & {[}0.34; 0.97895463741052202{]} & {[}0.35; 0.97901327197545851{]} &
{[}0.36; 0.97907199917004073{]} \\ \midrule
 {[}0.37; 0.97913480206897852{]} & {[}0.38; 0.97920757809296877{]} & {[}0.39; 0.97925580726117856{]} & {[}0.40; 0.97932396724121917{]} \\ \midrule
{[}0.41; 0.97938342003468337{]} & {[}0.42; 0.97944120549454572{]} & {[}0.43; 0.97948467967434194{]} & {[}0.44; 0.97954373107267778{]} \\ \midrule
 {[}0.45; 0.97959954043346009{]} &
{[}0.46; 0.97964355495284927{]} & {[}0.47; 0.97968756947223590{]} & {[}0.48; 0.97973988978307669{]} \\ \midrule
 {[}0.49; 0.97977845459171642{]} & {[}0.50; 0.97983410956975669{]} &
{[}0.51; 0.97987651850864055{]} & {[}0.52; 0.97992034776874026{]} \\ \midrule
 {[}0.53; 0.97998287277874474{]} & {[}0.54; 0.98001592612353894{]} & {[}0.55; 0.98006724294656333{]} &
{[}0.56; 0.98011216832412273{]} \\ \midrule
 {[}0.57; 0.98013869127898290{]} & {[}0.58; 0.98018823270049082{]} & {[}0.59; 0.98021707139646241{]} & {[}0.60; 0.98024144843122751{]} \\ \midrule
{[}0.61; 0.98026273781117690{]} & {[}0.62; 0.98028804114238566{]} & {[}0.63; 0.98031082803491987{]} & {[}0.64; 0.98034139581758906{]} \\ \midrule
 {[}0.65; 0.98037344567456941{]} &
{[}0.66; 0.98037471161304368{]} & {[}0.67; 0.98039501294345210{]} & {[}0.68; 0.98042110362663826{]} \\ \midrule
 {[}0.69; 0.98044646871094343{]} & {[}0.70; 0.98045438854554323{]} &
{[}0.71; 0.98046409921993705{]} & {[}0.72; 0.98048478650719761{]} \\ \midrule
 {[}0.73; 0.98048520334059863{]} & {[}0.74; 0.98049466700260646{]} & {[}0.75; 0.98049752308331006{]} &
{[}0.76; 0.98050144440492548{]} \\ \midrule
 {[}0.77; 0.98050379102258456{]} & {[}0.78; 0.98050913266541517{]} & {[}0.79; 0.98050024021954696{]} & {[}0.80; 0.98050201562106531{]} \\ \midrule
{[}0.81; 0.98049942199102114{]} & {[}0.82; 0.98049304598382936{]} & {[}0.83; 0.98048345881562826{]} & {[}0.84; 0.98046934823312271{]} \\ \midrule
 {[}0.85; 0.98045710568178235{]} &
{[}0.86; 0.98045168684758111{]} & {[}0.87; 0.98042429934937170{]} & {[}0.88; 0.98042289446643149{]} \\ \midrule
 {[}0.89; 0.98040929334697202{]} & {[}0.90; 0.98039153933178547{]} &
{[}0.91; 0.98037545265019832{]} & {[}0.92; 0.98036875243925059{]} \\ \midrule
 {[}0.93; 0.98033790676764754{]} & {[}0.94; 0.98032297795661794{]} & {[}0.95; 0.98029230210603036{]} &
{[}0.96; 0.98025315064297591{]} \\  \midrule
{[}0.97; 0.98021506442083295{]} & {[}0.98; 0.98019269436169787{]} & {[}0.99; 0.98015337307762984{]} & {[}1.00; 0.98012437999891566{]} \\
\bottomrule
\end{tabular}
\end{table}
 In each cell of the previous table, the $\delta$  and the corresponding S-index value appear inside the square brackets. \\

As already mentioned in Section \ref{sec5}, the Boat image was one of those  examined in  \cite{Mesiar}  which is the paper  that originated the comparison. \\
If we examine the shapes of the graphs of the \cite[Figure 12]{Mesiar} and that of Figure \ref{delta-fig}.4 we can observe that the qualitative curves are analogous.
 The maximum in the present paper is obtained for a larger value of the parameter $\delta$, but this may depend both on the
floating-point number format
  and on the fact that in this study we assume that the pixels outside the image have costant 
 value equal to zero, in fact they do not provide additional information (no boundary conditions),  while in \cite{Mesiar} this is not specified.
 In any case the difference of the S-index, in Table \ref{boat-num}, in the interval  between the ``$\delta_S^{\max}$'' of the
 two papers is less than $3 \cdot 10^{-4}$ thus we can conclude that the results obtained here confirm those in \cite{Mesiar}.\\

\clearpage
\subsection{magnification 5}
\phantom{a}\\
In this subsection the tables of the two similarity indices are presented  for the baboon and the mountain images; as   we said before, due to the size of the images, only these two could be taken into consideration for the magnification $R=5$. In this case the values  for which the maximum of the similarity indices is reached
are not   highlighted.
\begin{table}[h!]
\small
\caption{ \footnotesize  The numerical values of the PSNR-index for the baboon image R=5}\label{bsk-psnr-5}
\begin{tabular}{rccccc}
\toprule
N  & w=5                                      & w=10                                     & w=15                                     & w=20                                     & w=25                                     \\ \toprule
2  & 20,7329147502973                        & 21,7161082496651                     & 21,979365824288           &  22,0558557949496 & 22,0775418477414 \\ \midrule
3  & 20,3467942537683                         & 21,4828541953169                      & 21,8662212638324          & 22,0074149690313      & 22,0615063203985                         \\ \midrule
4  & 20,0415812414765                         & 21,2811780337                            & 21,7516133283618            & 21,9493699991135        & 22,021982797865                          \\ \midrule
5  & 19,8011752130362                         & 21,1064969190485                         & 21,6437660594206          & 21,8788040512415                         & 21,9870684510821                         \\ \midrule 
6  & 19,5936090626771                         & 20,9560372953294                         & 21,5387139742664       & 21,8223810785381       & 21,9529653869713                         \\ \midrule
7  & 19,4110677834914                         & 20,8227533041516                         & 21,4372167850169             & 21,7502537170246                         & 21,9079647824159                         \\ \midrule
8  & 19,2550258194904                         & 20,7062361523255                         & 21,3526599132005          & 21,6877191556461                         & 21,8689832613072                         \\ \midrule
9  & 19,1092740466907                         & 20,5972670626762                         & 21,2666457824789         & 21,6286366213533        & 21,8290361881159                         \\ \midrule
10 & 18,974699199968                          & 20,4979200147237                         & 21,1894189563848        & 21,5703646320501           & 21,7851824151841                         \\ \midrule
11 & 18,8531097177987                         & 20,403559417055                          & 21,113408433493        & 21,5085279324257        & 21,7433900697263                         \\ \midrule
12 & 18,7426431897043 & 20,3176158516323 & 21,0436294630194 & 21,4509224927457 & 21,7015964927653 \\
\bottomrule\\
\end{tabular}
\end{table}

\begin{table}[h!]
\footnotesize
\caption{ \footnotesize  The numerical values of the S-index for the baboon image R=5}\label{bsk-sindex-5}
 \begin{tabular}{r c c c c c }
\toprule
N  & w=5                                       & w=10                                      & w=15                                      & w=20                                      & w=25                                      \\
\toprule
2  & 0,931130741569984 & 0,939363412262253 & 0,941239854957746 & 0,941674318324023 & 0,941778411018386 \\ \midrule
3  & 0,92756891391697                          & 0,937563742451998                         & 0,940480870856609                         & 0,941408538194208                         & 0,941715991586946                         \\ \midrule
4  & 0,924526946649478                         & 0,935899284588883                         & 0,939632268132166                         & 0,941022804200496                         & 0,941476867871331                         \\ \midrule
5  & 0,922040482167492                         & 0,934411108849537                         & 0,938786801456453                         & 0,940544194917489                         & 0,941288343095793                         \\ \midrule
6  & 0,919811897384866                         & 0,933090651408584                         & 0,937969317984786                         & 0,940145615185713                         & 0,941051511108095                         \\ \midrule 
7  & 0,917812666319892                         & 0,931922835108668                         & 0,937176983211585                         & 0,939616527579889                         & 0,94074592728287                          \\ \midrule
8  & 0,916027108728919                         & 0,930876540697016                         & 0,936484760763205                         & 0,939115604104002                         & 0,940463622588597                         \\ \midrule
9  & 0,914360238520628                         & 0,929858832575705                         & 0,935755026347333                         & 0,938689704563101                         & 0,940181438511583                         \\ \midrule
10 & 0,912790043045284                         & 0,928941839865511                         & 0,935106346729387                         & 0,938211517440502                         & 0,939881885549299                         \\ \midrule
11 & 0,91133539890389                          & 0,928051202026371                         & 0,934448982668807                         & 0,937729832417396                         & 0,939541986113938                         \\ \midrule
12 & 0,910025374855824 & 0,927229798493791 & 0,933838599030538 & 0,937278603252143 & 0,939234834264347 \\  \bottomrule     \\                               
\end{tabular}
\end{table}
}
\clearpage


\begin{table}[h]
\small
\caption{ \footnotesize  The numerical values of the PSNR-index for the mountain image R=5}\label{msk-psnr-5}
\begin{tabular}{rcccc}
\toprule
$N$  & $ w=5 $             &         $w=10$         &       $w=15$                                     & $w=20$             \\ \toprule
2  & 19,4878493801042 & 20,5194643201115 & 20,3445774937346                         & 20,0352078339142 \\ \midrule
3  & 18,8959966960253 & 20,4106079549201 & 20,5078456050646                         & 20,2861784442023 \\  \midrule
4  & 18,4349923092873 & 20,2165121177552 & 20,5479617447969                         & 20,4262893816383 \\ \midrule
5  & 18,0690123314328 & 20,0065933263976 & 20,5227982717099                         & 20,5094371835843 \\  \midrule
6  & 17,7658779676526 & 19,8011176515348 & 20,4630175378557                         & 20,54521178183   \\  \midrule
7  & 17,5141848124694 & 19,6066833837114 & 20,3843693694262                         & 20,557343052487  \\  \midrule
8  & 17,2973630657787 & 19,4295932034864 & 20,2978171756184                         & 20,5473046959311 \\  \midrule
9  & 17,1092673976789 & 19,2613049769722 & 20,2016656072662                         & 20,5215375241213 \\  \midrule
10 & 16,9415207481794 & 19,1093674672082 & 20,1054519195362                         & 20,4886717646383 \\  \midrule
11 & 16,7911714644458 & 18,9661225386627 & 20,0125175550292                         & 20,4462480907923 \\  \midrule
12 & 16,655095783726  & 18,8313407159141 & 19,9169980743677                         & 20,402755638215 \\ \bottomrule
\\ \\
\end{tabular}
\end{table}


\begin{table}[h]
\small
\caption{ \footnotesize  The numerical values of the S-index for the mountain image R=5}\label{msk-sindex-5}
\begin{tabular}{rcccc}
\toprule
$N$  & $ w=5 $             & $w=10$                    & $w=15$                    & $w=20$             \\ \toprule
2  & 0,926103965141612 & 0,937437124183007 & 0,937670007262164 & 0,935860958605665 \\ \midrule
3  & 0,919899927378358 & 0,935601263616558 & 0,938177777777779 & 0,937387596223674 \\ \midrule
4  & 0,914731169208424 & 0,933381859114015 & 0,937776441539579 & 0,938012273057373 \\ \midrule
5  & 0,910410399419027 & 0,931204371822803 & 0,937041626724764 & 0,938218997821352 \\ \midrule
6  & 0,906666129266522 & 0,929122004357298 & 0,936115889615106 & 0,938098634713145 \\ \midrule
7  & 0,903401045751634 & 0,927175061728395 & 0,935151154684096 & 0,93781179375454  \\ \midrule
8  & 0,900529339143064 & 0,925371082062454 & 0,934170806100217 & 0,937403877995643 \\ \midrule
9  & 0,8979444734931    & 0,9236483805374      & 0,933156165577341   & 0,936923384168482 \\ \midrule
10 & 0,895572127814088 & 0,922042527233116 & 0,932153972403777 & 0,936412549019608 \\ \midrule
11 & 0,893384836601306 & 0,920517618010166 & 0,931215003631082 & 0,935867494553378 \\ \midrule
12 & 0,891373623819897 & 0,919074669571532 & 0,930259288307916 & 0,935333710965867 \\
\bottomrule                  
\end{tabular}
\end{table}

\section{Conclusions} 
 
In this article we compared a construction method
of an interval-valued fuzzy set starting from fuzzy sets, introduced in \cite{Mesiar} with the SK algorithm 
and the well-known bicubic method for digital image processing. 
These algorithms  were  compared with the use of the PSNR and the likelihood S indices, 
as well as, by  analysing
 the corresponding processing CPU time. \\
From the numerical results provided in Section \ref{sec5} it seems to be clear that:
\begin{itemize}
\item Based on the analysis of Tables \ref{tab1} and  \ref{tab2}, it seems that the maximum values of the PSNR and
 likelihood index S are both substantially better in the case of the application of the SK method with sufficiently high $w$,
 with respect to 
 other two considered methods. Only when the scaling factor is equal to 3, 
and we consider the "boat",  does  the fuzzy-algorithm seem to provide better reconstruction results,
  at least for the PSNR index. 
The same consideration can also  be applied when the scaling factor is equal to 5. \\
 The fuzzy-type algorithm seems to perform substantially better than the bicubic method.
\item The CPU analysis given in Tab. \ref{tab3},
 performed only for the best approximations,
 shows that the  bicubic method has the  most rapid execution, the mean CPU time employed by the fuzzy-type algorithm
 is reasonable in term of applicability of the method, while, as we already known, the CPU time is the weak point of the SK algorithm. 
 The higher CPU time 
seems to be the price to pay  to  obtain more accurate results.

\end{itemize}

{\bf \noindent Author's contribution}
{\small All  authors    contributed  equally  to  this  work  for  writing,  reviewing  and  editing. All authors have read and agreed
 to the published version of the manuscript.}\\

{\bf \noindent Conflict of interest} {\small  The authors declare no conflicts of interest.}\\

{\bf \noindent Copyright }  
{\small
The figures  (baboon, boat, mountain) are  contained in the repository 
https://links.uwaterloo. ca/Reposi\-to\-ry.html and they belong to the Grayscale Set 2 (The Waterloo Fractal Coding and Analysis Group).
This set of images was formally part of the BragZone repository  https://links.uwaterloo.ca/old website/ bragzone.base.html 
(this resource is intended for researchers and graduate students), \cite{repository}.
The last image (city) was contained in the Data Set given in the article \cite{foto}, by  M. Castro, DM. Ballesteros, D. Renza, 
 under license CC BY 4.0.\\}


{\bf \noindent Data Availability Statement:} 
{\small
All the data generated for this study were stored in our laboratory and are
not publicly available. Researchers who wish to access the data directly  contacted the corresponding
author.}\\

{\bf \noindent Funding}
{\small This research has been accomplished within the 
UMI Group TAA- “Approximation Theory and Applications”, 
the group RITA - "Research ITalian network on Approximation",
 the G.N.AM.P.A. group of   INDAM and the University of Perugia.
This study was partly funded by:\\
-   "National Innovation Ecosystem grant ECS00000041 - VITALITY", (European Union - NextGenerationEU) under  MUR;\\
- Research project of MIUR (Italian Ministry of Education, University and Research)
 Prin 2022  “Nonlinear differential problems with applications to real phenomena” (Grant Number: 2022ZXZTN2, CUP J53D23003920006 );\\
-PRIN 2022 PNRR: “RETINA:
REmote sensing daTa INversion with multivariate functional modelling for essential climAte variables
characterization” funded by the European Union under the Italian National Recovery and
Resilience Plan (NRRP) of NextGenerationEU, under the MUR (Project Code: P20229SH29,
CUP: J53D23015950001); \\
-Gnampa Project 2023 "Approssimazione co\-strut\-tiva e astratta mediante operatori di tipo sampling e loro applicazioni".
}
%


\end{document}